\documentclass[12pt,twoside]{amsart}
\usepackage{amssymb,amsmath,amsfonts,amsthm}
\usepackage{mathrsfs}
\usepackage[X2,T1]{fontenc}
\usepackage[applemac]{inputenc}
\usepackage{cite}
\usepackage{latexsym}
\usepackage[dvips]{graphicx}
\usepackage{comment}
\usepackage{ulem}
\def\Im{\mathop{\rm Im}\nolimits}

\def\Im{\mathop{\rm Im}\nolimits}

\def\R{\mathbb R}
\def\C{\mathbb C}
\def\N{\mathbb N}

\def\ds{\displaystyle}

\newcommand\dslash{d\llap {\raisebox{.9ex}{$\scriptstyle-\!$}}}
\usepackage{colortbl}

\newcommand{\beqsn}{\arraycolsep1.5pt\begin{eqnarray*}}
	\newcommand{\eeqsn}{\end{eqnarray*}\arraycolsep5pt}
\newcommand{\beqs}{\arraycolsep1.5pt\begin{eqnarray}}
	\newcommand{\eeqs}{\end{eqnarray}\arraycolsep5pt}

\newtheorem{Th}{Theorem}[section]
\newtheorem{Rem}[Th]{Remark}
\newtheorem{Ex}[Th]{Example}
\newtheorem{Lemma}[Th]{Lemma}
\newtheorem{Def}[Th]{Definition}
\newtheorem{Prop}[Th]{Proposition}
\newtheorem{Cor}[Th]{Corollary}

\catcode`\@=11

\renewcommand{\section}%
{\setcounter{equation}{0}\@startsection {section}{1}{\z@}{-3.5ex plus -1ex
		minus -.2ex}{2.3ex plus .2ex}{\Large\bf}}

\title{Schr\"odinger type equations with singular coefficients and lower order terms}
\author[Arias Junior]{Alexandre Arias Junior}
\address{Alexandre Arias Junior\\
	Dipartimento di Matematica ``G. Peano" \\Universit\`a di Torino\\
	Via Carlo Alberto 10\\
	10123 Torino\\
	Italy}
	\email{alexandre.ariasjunior@unito.it}
\author[Ascanelli]{Alessia Ascanelli}
\address{Alessia Ascanelli\\
	Dipartimento di Matematica ed Informatica\\Universit\`a di Ferrara\\
	Via Machiavelli 30\\
	44121 Ferrara\\
	Italy}
\email{alessia.ascanelli@unife.it}
\author[Cappiello]{Marco Cappiello}
\address{Marco Cappiello\\
	Dipartimento di Matematica ``G. Peano" \\Universit\`a di Torino\\
	Via Carlo Alberto 10\\
	10123 Torino\\
	Italy}
\email{marco.cappiello@unito.it}
\author[Garetto]{Claudia Garetto}
\address{Claudia Garetto\\
	School of Mathematical Sciences\\Queen Mary University of London\\
Mile End Road, London E1 4NS\\ UK}
\email{c.garetto@qmul.ac.uk}

\thanks{The last author was supported by the
EPSRC grant EP/V005529/2.}
\thanks{On behalf of all authors, the corresponding author states that there is no conflict of interest. This manuscript has no associated data.}

\begin{document}

\begin{abstract}
	In this paper we investigate the well-posedness of the Cauchy problem for a Schr\"odinger operator with singular lower order terms. We allow distributional coefficients and we approach this problem via the regularising methods at the core of the theory of very weak solutions. We prove that a very weak solution exists and it is unique modulo negligible perturbations. Very weak solutions converge to classical solutions when the equation coefficients are regular enough. \end{abstract}
\maketitle
\section{Introduction}
The problem of existence of a unique solution to the initial value problem for Schr\"odinger type equations in suitable spaces of functions or distributions has been widely investigated in literature in the past years: the Schr\"odinger equation is indeed the fundamental equation in Quantum Mechanics. Many papers deal with it, in the linear and nonlinear case, in the deterministic and stochastic case, for flat or asymptotically flat metrics, with or without potentials. Several results of existence and uniqueness of a (possibly weak) solution have been given under {\sl{suitable assumptions of regularity on the coefficients}} of the Schr\"odinger operator. The aim of this paper is to go beyond such assumptions of regularity. Although the literature on this topic is huge, let us try to summarise some basic facts.
Consider the initial value problem
\begin{equation}\label{CP}
	\begin{cases}
		S u(t,x) = f(t,x), \quad\, t \in [0,T],\, x \in \R^{n}, \\
		u(0,x) = g(x), \quad \quad \,\, x \in \R^{n},
	\end{cases}
\end{equation}
associated to a class of Schr\"odinger type operators of the form
$$
S = D_t - a(t) \triangle_x + \sum_{j=1}^{n}c_j(t,x)D_{x_j} + c_0(t,x), \quad t \in [0,T],\, x \in \R^{n},
$$
where $\triangle _x:=\sum_{j=1}^{n}\partial^2_{x_j}$ and $D:=-i\partial.$ The well-posedness of the problem \eqref{CP} is well understood in the case when the coefficients are at least continuous with respect to $t$ and of class $\mathcal{B}^{\infty}(\R^n)$ with respect to $x$, where $\mathcal{B}^{\infty}(\R^n)$ stands for the space of all complex-valued smooth and uniformly bounded functions together with all their derivatives. When all the coefficients $c_j$ are real-valued, the lower order terms define a self-adjoint operator and, in this case, well-posedness is obtained in a straightforward way. On the contrary, i.e., when some of the coefficients $c_j$ have a non identically zero imaginary part, the situation is more delicate. We summarise some of the known results:

\begin{itemize}
	\item When the coefficients do not depend on the time variable and $a(t) = 1$, a necessary condition to $H^\infty= \bigcap_{m \in \R} H^m(\R)$ well-posedness is that there exist constants $M,N >0$ such that
	$$
	\sup_{x \in \R^{n}, \omega \in S^{n-1}} \left| \sum_{j=1}^{n} \int_{0}^{\rho} Im\, c_j(x+\omega\theta) d\theta \right| \leq M\log(1+\rho) +N, \quad \forall\, \rho \geq 0, 
	$$
	see \cite{ichinose_remarks_cauchy_problem_schrodinger_necessary_condition};
	
	\item the previous condition is also sufficient in the case $n=1$, otherwise technical assumptions on the derivatives of $c_j$ have to be added, see \cite{I2};
	
	\item if $a(t)$ is real-valued and continuous, $c_j, c_0$ are continuous in $t$ and of class $\mathcal B^\infty(\R^n_x)$ and $|Im\, (c_j(t,x))|\leq C\langle x\rangle^{-\sigma}$ for every $x\in\R^n$ and $t \in [0,T]$, where $\langle x\rangle:=(1+|x|^2)^{1/2}$, then the Cauchy problem \eqref{CP} is well posed in $H^m(\R^n)$ for every $m\in\R$ if $\sigma>1$, in $H^\infty(\R^n)$ if $\sigma=1$ and in Gevrey classes of index $s\leq1/(1-\sigma)$ if $0<\sigma<1$, see \cite{KB};
	
	\item the case when $a(t)$ vanishes of some order has been investigated in \cite{CicRei}, where the authors obtained a well posed Cauchy problem either in $H^\infty$ or in Gevrey classes depending on the vanishing order of $a(t)$ and on the corresponding Levi conditions that they ask on the lower order terms $c_j$;
	\item the Gevrey well posedness under assumption $|Im\, c_j(t,x)|\leq C\langle x\rangle^{-\sigma}$, $0<\sigma<1$, in the case of a Gevrey index $s> 1/(1-\sigma)$ has been investigated in \cite{ACR}, under exponential decay assumptions on the data;
	
	\item the Cauchy problem for Schr\"odinger-type equations in Gelfand-Shilov spaces has been studied in \cite{Arias_GS} and \cite{scncpp5}.
	
\end{itemize}

Concerning Schr\"odinger equations with singular  coefficients, there are several works treating the case of singular potentials, that is when $c$ belongs to $L^p$ or modulation spaces, see for instance \cite{CardonaRuzhansky, CorderoNicolaRodino, DanconaPierfeliceVisciglia}. In these papers, the coefficients of the lower order terms are real-valued. \\
The aim of this paper is to analyse a very singular situation, that is when all the coefficients of the lower order terms may be distributions.
Namely we focus on operators of the form 
\begin{equation}\label{true_operator_S}
S = D_t - a(t) \triangle_x + \sum_{j=1}^{n}a_j(t)b_j(x)D_{x_j} + a_0(t)b_0(x), \quad t \in [0,T],\, x \in \R^{n},
\end{equation}
where $a(t)$ is a real-valued bounded function which never vanishes, $a_j$ and $a_0$ are compactly supported distributions on an interval $I$ containing $[0,T]$ and $b_{j}$, $b_0$ are certain tempered distributions that we shall define later on. The choice of coefficients defined as tensor products is motivated by their different behavior with respect to time and space variables. 
\\
The situation described above is clearly very general and therefore requires a redefinition of the notion of solution of the problem \eqref{CP} because of intrinsic problems due to multiplication of distributions. In the recent papers \cite{G,GR}, the last author et al. treated  hyperbolic equations with discontinuous functions or compactly supported distributions as coefficients and introduced the concept of {\sl{very weak solution}}.  Note that evolution equations with non-regular coefficients appear frequently in geophysics when modelling wave transmission through the  Earth subsoil which has a multilayered and therefore discontinuous structure. However, allowing distributional coefficients leads to the major problem of identifying a reasonable notion of solution since the equation operator might fail to be well-defined when the coefficients are less than continuous. Following the approach introduced in \cite{GR} we will replace the operator $S$ with a family of regularised operators $S_\varepsilon$, $\varepsilon\in(0,1]$, obtained via convolving the irregular coefficients with a net of mollifiers $\varphi_{\omega(\varepsilon)}(\cdot)=\omega(\varepsilon)^{-n}\varphi(\cdot/\omega(\varepsilon))$, where $\omega(\varepsilon)$ is a positive scale converging to $0$ as defined later in the paper. We therefore look at the net of solutions $(u_\varepsilon)_\varepsilon$ of the regularised problem and we provide a qualitative analysis of $(u_\varepsilon)_\varepsilon$ with respect to the parameter $\varepsilon$ by analysing its limiting behaviour as $\epsilon$ tends to $0$. In a nutshell, this means to find a {\sl{very weak solution}} for our Cauchy problem. For recent applications of the theory of very weak solutions and for some new insights provided via numerical experiments we refer the reader to  \cite{DGL, DGL2, MRT19, GS24}.  The paper is organised as follows. In Section 2 we describe our main result concerning the existence of a very weak solution for the Cauchy problem \eqref{CP} and we provide the needed preliminaries. The construction and analysis of the regularised problem is given in Section 3. Few results concerning pseudo-differential operators employed in the paper can be found in Section 4. The proof of our main result is spread throughout Sections 5, 6 and 7, with a different approach for the 1-dimensional case and $n$-dimensional case. The paper ends with Section 8, where we discuss the uniqueness of the very weak solution and we prove consistency with the classical theory in case of regular coefficients.

%
%
%
%


\section{Main result}\label{Main result}
In this section we state our main result concerning the existence of a very weak solution for the Cauchy problem \eqref{CP}. The statement requires some preliminaries. \\

Let $\varphi$ be a Schwartz function such that $\int_{\R^{n}} \varphi = 1$. Given a positive scale $\omega(\varepsilon)$, $\varepsilon \in (0,1]$, i.e. $\omega$ is positive and bounded, $\omega(\varepsilon) \to 0$ as $\varepsilon \to 0^{+}$ and $\omega(\varepsilon) \geq c \varepsilon^{r}$, for some $c, r > 0$, we define
$$
\varphi_{\omega(\varepsilon)}(x) = \frac{1}{(\omega(\varepsilon))^{n}} \varphi \left( \frac{x}{\omega(\varepsilon)} \right).
$$
As usual, tempered distributions can be regularised via convolutions with the net $(\varphi_{\omega(\varepsilon)})_{\varepsilon}$. The following result, which deals with regularisations of different kind of distributions and functions, follows from standard arguments, as in \cite{GR, GKOS01, Obe92},  and for this reason we shall omit the proof. 

\begin{Prop}\label{the_evil_that_men_do}\leavevmode
	\begin{itemize}
		\item If $u \in \mathscr{S}'(\R^{n})$ then for any $\beta \in \N_0^n$ there exists $N_1(n,u) = N_1$ and $N_2(u) = N_2$ such that 
		$$
		|\partial^{\beta}_{x} (\varphi_{\omega(\varepsilon)} \ast u)(x)| \leq C \omega(\varepsilon)^{-N_1-|\beta|} \langle x \rangle^{N_2};
		$$
		
		\item If $u \in \mathcal{E}'(\R^{n})$ then for any $\beta \in \N_0$ there exists $N(n,u)$ such that 
		$$
		|\partial^{\beta}_{x} (\varphi_{\omega(\varepsilon)} \ast u)(x)| \leq C \omega(\varepsilon)^{-N-|\beta|};
		$$
		
		\item If $u \in \mathcal{B}^{\infty}(\R^{n})$ then for any $\beta \in \N_0^n$ there exists $c > 0$ such that 
		$$
		|\partial^{\beta}_{x} (\varphi_{\omega(\varepsilon)} \ast u)(x)| \leq c;
		$$
		
		\item  If $u \in \mathcal{B}^{\infty}(\R^{n})$ and $\varphi$ has all moments vanishing, i.e. $\int x^{\alpha} \varphi(x) dx = 0$ for all $\alpha \neq 0$, then for any $\beta \in \N_0^n$ and any $q \in \N_0$ there exists $c > 0$ such that 
		$$
		|\partial^{\beta}_{x} (\varphi_{\omega(\varepsilon)} \ast u - u)(x)| \leq c (\omega(\varepsilon))^{q}.
		$$
	\end{itemize}
\end{Prop}

We point out that from $\omega(\varepsilon)\ge c\varepsilon^r$ it follows that $\omega(\varepsilon)^{-N_1-|\beta|} $ and $\omega(\varepsilon)^{-N-|\beta|}$ in the proposition above can be replaced by $\varepsilon^{-M}$ for some $M$ depending on $u, n, \beta$ and the scale $\omega$.

In the sequel we consider the concept of $H^{\infty}$-moderateness and $H^{\infty}$-negligibility of nets of $H^{\infty}(\R^{n})$ functions and analogously the same concepts with $H^\infty$ replaced by $\mathcal{B}^\infty$. Note that we will work with nets of functions in the variables $t\in[0,T]$ and $x\in\R^n$ and that we will always assume boundedness with respect to $t\in[0,T]$. $H^\infty$- or $\mathcal{B}^\infty$-estimates will be therefore considered with respect to the space variable $x$.

\begin{Def} \label{defmoderateness}
	\leavevmode
	\begin{itemize}
		\item[(i)] Let $(v_{\varepsilon})_{\varepsilon} \in \{C([0,T];H^{\infty}(\R^{n}))\}^{(0,1]}$. We say that the net $(v_\varepsilon)_\varepsilon$ is $H^{\infty}$-moderate if for any $m\in\N_0$ there exists $N \in \N_0$ and $C > 0$ such that 
		$$
		\|v_{\varepsilon}(t,\cdot)\|_{H^m} \leq C \varepsilon^{-N}, \quad \forall t\in [0,T], \epsilon\in (0,1].
		$$
		 \item[(ii)] Let  $(v_{\varepsilon})_{\varepsilon} \in \{C([0,T];H^{\infty}(\R^n))\}^{(0,1]}$ We say that the net $(v_{\varepsilon})_{\varepsilon}$ is $H^{\infty}$-negligible if for any $m\in\N_0$ and for any $q\in\N_0$ there exists $C > 0$ such that 
		$$
		\|v_{\varepsilon}(t,\cdot)\|_{H^m} \leq C \varepsilon^{q}, \quad \forall t\in [0,T], \epsilon\in (0,1].
		$$
		\item[(iii)] Let $(v_{\varepsilon})_{\varepsilon} \in \{C([0,T];\mathcal{B}^{\infty}(\R^{n}))\}^{(0,1]}$. We say that the net $(v_\varepsilon)_\varepsilon$ is $\mathcal{B}^{\infty}$-moderate if for all $\beta \in \N_0$ there exist $N \in \N_0$ and $C > 0$ such that 
		$$
		\sup_{x\in\R^{n}}|\partial^{\beta}_{x}v_{\varepsilon}(t,x)| \leq C \varepsilon^{-N},\quad \forall t\in [0,T], \epsilon\in (0,1].
		$$
		\item[(iv)] Let $(v_{\varepsilon})_{\varepsilon} \in \{C([0,T];\mathcal{B}^{\infty}(\R^{n}))\}^{(0,1]}$. We say that the net $(v_\varepsilon)_\varepsilon$ is $\mathcal{B}^{\infty}$-negligible if for all $\beta \in \N_0$ and for all $q\in\N_0$  there exists $C > 0$ such that 
		$$
		\sup_{x\in\R^{n}}|\partial^{\beta}_{x}v_{\varepsilon}(t,x)| \leq C \varepsilon^{q},\quad \forall t\in [0,T], \epsilon\in (0,1].
		$$
	\end{itemize}
\end{Def}

\begin{Rem}
	The above definitions $(iii)$ and $(iv)$ are a special case of $C^{\infty}$-moderateness and $C^\infty$-negligibility, respectively where we ask a uniform estimate on the whole $\R^n$ instead of on compact sets. See \cite{GKOS01, Obe92}. Note also that by Sobolev's embedding theorem $H^\infty(\R^n)\subset \mathcal{B}^\infty(\R^n)$. Hence, $(i)$ implies $(iii)$ and $(ii)$ implies $(iv)$.
\end{Rem}


Now we introduce the concept of very weak solution that we are interested in. 

\begin{Def}
	\label{def_vw}
	The net $(u_\varepsilon)_{\varepsilon} \in \{C([0,T];H^{\infty}(\R^{n}))\}^{(0,1]}$ is a $H^{\infty}$ very weak solution for the Cauchy problem \eqref{CP} if  there exist 
	\begin{itemize}
		\item $\mathcal{B}^{\infty}$-moderate regularisations $(a_{\varepsilon})_{\varepsilon}$ and $(a_{j,\varepsilon})_{\varepsilon}$ of $a$ and $a_j$, $j = 0, 1, \ldots, n$,
		
		\item $\mathcal{B}^{\infty}$-moderate regularisations $(b_{j,\varepsilon})_{\varepsilon}$ of $b_j$, $j = 0, 1, \ldots, n$,
		
		\item $H^{\infty}$-moderate regularisations $(f_{\varepsilon})_{\varepsilon}, (g_{\varepsilon})_{\varepsilon}$ of the Cauchy data $f$ and $g$,
	\end{itemize}
	such that, for every fixed $\varepsilon$, $u_{\varepsilon}$ solves the Cauchy problem
	\begin{equation}\label{regularized_CP}
		\begin{cases}
			S_{\varepsilon} v(t,x) = f_\epsilon(t,x), \quad t \in [0,T],\, x \in \R^{n}, \\
			v(0,x) = g_\epsilon(x), \quad \quad \,\, \, x \in \R^{n},
		\end{cases}
	\end{equation}
for the regularised operator
\begin{equation}\label{regoperator}
	S_\varepsilon = D_t - a_\varepsilon(t) \triangle_x + \sum_{j=1}^{n}a_{j,\varepsilon}(t)b_{j,\varepsilon}(x)D_{x_j} + a_{0,\varepsilon}(t)b_{0,\varepsilon}(x), \quad t \in [0,T],\, x \in \R^{n},
\end{equation}
and $(u_\varepsilon)_\varepsilon$ is $H^{\infty}$-moderate.
\end{Def}
\begin{Rem}
The definition above is not affected by replacing the parameter interval $(0,1]$ with any smaller interval $(0,\varepsilon_0]$. This is actually quite natural when regularising distributions on a bounded interval $[0, T]$ as we will see in the next section.
\end{Rem}

To obtain a very weak solution for \eqref{CP} we need to solve the regularised problem \eqref{regularized_CP}, \eqref{regoperator} in $H^{\infty}(\R^{n})$. So, in order to apply the techniques coming from classical theory, for instance the ones used in \cite{KB}, we need the following behavior for the regularised coefficients and regularised data:
\begin{itemize}
	\item $a_{\varepsilon}(t)$ should be a real-valued continuous function on $[0,T]$ which never vanishes; 
	
	\item $a_{j,\varepsilon}(t)$, $j = 0, 1, \ldots, n$, should be continuous functions on $[0,T]$;
	
	\item $b_{j,\varepsilon}(x)$, $j = 0, 1, \ldots, n$, should belong to $\mathcal{B}^{\infty}(\R^{n})$;
	
	\item $|Im\, (a_{j,\varepsilon}(t)b_{j,\varepsilon}(x))| \lesssim \langle x \rangle^{-\sigma}$, $j = 1, \ldots, n$, uniformly in $t$ and for some $\sigma \geq 1$;
	
	\item $f_\varepsilon \in C([0,T]; H^{\infty}(\R^{n}))$ and $g_{\varepsilon} \in H^{\infty}(\R^{n})$.
\end{itemize}
We therefore need a suitable set of hypotheses on the coefficients $a, a_j, b_j$ and on the data $f, g$ in such a way that once we regularise them, the obtained regularisations satisfy the conditions above and the needed moderateness assumptions.

Regarding the coefficients $a, a_j$, $j = 0, 1, \ldots, n$ we have the following quite natural hypotheses: 
\begin{equation*}
	a(t)\, \text{is a positive bounded function such that} \,\, 0 < C_{a} \leq a(t) \leq \tilde{C}_{a}; 
\end{equation*}
\begin{equation*}
	a_j \in \mathcal{E}'(I), j = 0, 1, \ldots, n, \,\, \text{where}\, I \, \text{is an open interval containing}\, [0,T].
\end{equation*}

In order to find suitable assumptions for the coefficients $b_j$, $j = 0, 1, \ldots, n$ and data $f, g$ we need to take into account both regularity and behavior at infinity. From Proposition \ref{the_evil_that_men_do} we know that regularisations of tempered distributions in general do not give functions in $\mathcal{B}^{\infty}(\R^{n})$. Therefore we are led to consider some subclass of tempered distributions. A first possibility based on classical distributional spaces could be assuming $b_j$ and $b_0$ in
$$
\mathcal{O}_{C}' = \mathcal{F}^{-1} \mathcal{O}_{M} = \{ u \in \mathscr{S}'(\R^{n}): \widehat{u} \in \mathcal{O}_{M} \},
$$
where
$$
\mathcal{O}_{M} = \left\{ f \in C^{\infty}(\R^{n}) : \text{for any}\,\, \alpha \in \N_0^n \,\, \text{there is}\, p \geq 0\,\, \text{such that} \,\, \sup_{x} \{\langle x \rangle^{-p} |\partial^{\alpha}_{x}f(x)|\} < \infty \right\}.
$$
In this case the regularised coefficients turn out to be Schwartz functions.

Since Schwartz regularity, in particular rapid decay, is much more than what we need, we introduce the following spaces: for any $i \in \N_0$ we consider
\begin{equation}
	H^{-\infty, i}(\R^{n}) = \left\{ u \in \mathscr{S}'(\R^{n}) :\,  \text{for all}\,\, |\beta|\leq i \,\, \text{there is} \,\, p_{\beta} \geq 0 \,\, \text{such that} \,\, \sup_{\xi \in \R^{n}} |\partial^{\beta}_{\xi} \widehat{u}(\xi)| \leq C_{\beta} \langle \xi \rangle^{p_{\beta}} \right\}.
\end{equation}
The motivation to consider the above space is the following, if $u \in H^{-\infty, i}(\R^{n})$ then $\widehat{u} \in C^{i}(\R^{n})$ and all the derivatives $\partial^{\beta}\widehat{u}$ up to order $i$ have at most a polynomial growth, hence the regularisation $u_{\varepsilon} = u \ast \varphi_{\varepsilon}$  will be a $\mathcal{B}^{\infty}(\R^{n})$ function with decay $\langle x \rangle^{-i}$. Indeed, take $\alpha \in \N_0^n$ and $|\beta| \leq i$. By straightforward computations we have
\begin{align*}
	x^{\beta} \partial^{\alpha}_{x} u_{\varepsilon}(x) = (-1)^{\beta} \frac{1}{\omega(\varepsilon)^{|\alpha|+n}}\int e^{ix\xi} D^{\beta}_{\xi} \{\widehat{u}(\xi) \widehat{\partial^{\alpha}\phi}(\omega(\varepsilon)\xi)\} \dslash\xi,
\end{align*}
where $\dslash \xi := \frac{d\xi}{(2\pi)^{n}}$. So, for all $\alpha \in \N^{n}_0$ and $|\beta| \leq i$ there exists $p = p(u, i) \geq 0$ such that 
\begin{align}\label{dec}
	|x^{\beta} \partial^{\alpha}_{x} u_{\varepsilon}(x)| \leq  C_{\alpha} \omega(\varepsilon)^{-n-|\alpha|-p},
\end{align}
which allow us to conclude our intuitive assertion. We also point out that if $i \geq n+1$, then $u_{\varepsilon}$ and all their derivatives belong to $L^{2}(\R^{n})$, so, in view of \eqref{dec} we conclude that $(u_{\varepsilon})_{\varepsilon}$ is a $H^{\infty}$-moderate net.  

\begin{Ex}
	For all $\xi \in \R$ we consider
	$$
	\rho(\xi) = \sin(e^{\xi^2})e^{-\xi^2}.
	$$
	We then have 
	$$
	\partial_{\xi} \rho(\xi) = \cos(e^{\xi^2}) 2\xi - \sin(e^{\xi^2})2\xi e^{-\xi^2}
	$$
	and any higher order derivative $\partial^{\alpha}_{\xi}\rho$ ($\alpha \geq 2$) cannot be bounded by a polynomial. Therefore $a = \mathcal{F}^{-1}(\rho) \in H^{-\infty, 1} (\R) - \mathcal{O}'_{C}.$
	
	Since $\rho$ has a super exponential decay at $|\xi| \to +\infty$, then $a \in C^{\infty}$. On the other hand, if we consider for instance 
	$$
	\lambda(\xi) = 
	\begin{cases}
		\sin(e^{\xi}-1)e^{-\xi} - \xi, \quad \xi \geq 0,\\
		0, \qquad \qquad \qquad \qquad \xi < 0,
	\end{cases}
	$$
	we have $b = \mathcal{F}^{-1}(\lambda) \in H^{-\infty, 1} (\R) - \mathcal{O}'_{C}$, but in this case $b$ it is not even a function. 
\end{Ex}

\begin{Rem}
For any $i = 0, 1, 2, \ldots$ the following inclusions hold:
	$$
	\mathcal{E}'(\R^{n}) \subset \mathcal{O}_{C}' \subset H^{-\infty,i+1}(\R^{n}) \subset H^{-\infty,i}(\R^{n}) \subset H^{-\infty}(\R^{n}) = \bigcup_{m \in \R} H^{m}(\R^{n}) \subset \mathscr{S}'(\R^{n}). 
	$$
\end{Rem}

We are finally ready to state the hypotheses that we are going to consider on the coefficients $b_j$, $j = 0, 1, \ldots, n$:
\begin{equation*}
	b_{j} \in H^{-\infty, 2}(\R^{n}),\, j =1, \ldots, n, \quad b_{0} \in H^{-\infty,0}(\R^{n}).
\end{equation*}
For the Cauchy data we shall ask the following:
\begin{equation*}
	f \in C([0,T], H^{-\infty, n+1}(\R^{n})), \, g\in H^{-\infty, n+1}(\R^{n}).
\end{equation*}
Note that we endow the space $H^{-\infty,n+1}(\R^{n})$ with the following notion of convergence (and related topology):  we say that $u_j \in H^{-\infty,n+1}(\R^{n})$ converges to $u \in H^{-\infty,n+1}(\R^{n})$ if there exist $p \geq 0$ and $C > 0$ such that 
$$
\sup_{\xi \in \R^{n}, |\beta| \leq n+1}|\partial^{\beta}_{\xi}\widehat{u}(\xi)| \langle \xi \rangle^{-p} \leq C, \sup_{\overset{\xi \in \R^{n}, |\beta| \leq n+1}{j \in \N_0}}|\partial^{\beta}_{\xi}\widehat{u}_j(\xi)| \langle \xi \rangle^{-p} \leq C
$$ 
and 
$$
\sup_{\xi \in \N_0, |\beta| \leq n+1} \langle \xi \rangle^{-p}|\partial^{\beta}_{\xi}(\widehat{u}_{j}-\widehat{u})(\xi)| \to 0 \,\, \text{as} \,\, j \to \infty.
$$
It easily follows from the arguments seen below that if $f \in C([0,T];H^{-\infty,n+1}(\R^{n}))$ then $f_{\varepsilon}(t,x) := \varphi_{\omega(\varepsilon)}(x) \ast f(t,x) \in C([0,T];H^{\infty}(\R^{n}))$ and there exists $p \geq 0$ such that 
$$
|\partial^{\beta}_{x} f_{\varepsilon}(t,x)| \leq C_{\beta} \varepsilon^{-n-p-|\beta|} \langle x \rangle^{-n-1}, \quad t \in [0,T], x \in \R^{n}.
$$

\smallskip


We are now ready to state our main result.

\begin{Th}\label{mainn}
Consider the Cauchy problem \eqref{CP} under the following hypotheses on the coefficients:
\begin{itemize}
	\item [(i)] $a(t)$ is a real-valued function satisfying 
	\begin{equation}
		0 < C_{a} \leq a(t) \leq \tilde{C}_{a}, \quad t \in [0,T],
	\end{equation}
	\item [(ii)] $a_j \in \mathcal{E}'(I)$, $j=0,1,\ldots,n$, where $I$ is an open interval containing $[0,T]$,
	\item [(iii)] $b_j \in H^{-\infty,2}(\R^{n})$, $j = 1, \ldots, n$, and $b_0 \in H^{-\infty,0}(\R^{n})$.
\end{itemize} 
Assume $f \in C([0,T], H^{-\infty, n+1}(\R^{n})), g\in H^{-\infty, n+1}(\R^{n})$. Then, \eqref{CP} admits a very weak solution of $H^\infty$ type.
\end{Th}

\begin{Rem}
	A decay like $\langle x \rangle^{-1}$ is enough to get a  solution in $H^{\infty}(\R^{n})$ for the regularised problem \eqref{regularized_CP}, see \cite{KB}, so in principle it would be sufficient to take the coefficients $b_j \in H^{-\infty,1}(\R^n)$. However, in this case the solution of the regularised problem would exhibit a loss of Sobolev regularity with respect to the initial data. This loss depends in general on the first order coefficients, and hence also on the parameter $\varepsilon$. This would make the dependence of the constants appearing in the energy estimates very difficult to control. For this reason we decided to take the coefficients $b_j$ in $H^{-\infty,2}(\R^{n})$ because under this assumption the solution of the regularised problem has the same regularity of the initial data. 
\end{Rem}

The proof of Theorem \ref{mainn} consists of several steps and it is organised as follows.
In Section \ref{rp} we regularise the Cauchy problem \eqref{CP}, constructing then a family of regularised associated Cauchy problems. These regularised problems will admit a unique classical solution in $H^\infty(\R^{n})$.
By solving these Cauchy problems in the $H^\infty$ framework, we will obtain a net of solutions $(u_{\varepsilon})_{\varepsilon \in (0,\varepsilon_{0}]}$ for some sufficiently small $\varepsilon_0>0$. Moreover, we will also derive energy estimates for the solutions $u_{\varepsilon}$ writing explicitly how the constants depend on the parameter $\varepsilon$. Since the regularised problem can be treated in a simpler way in space dimension $1$ we dedicate Section \ref{n=1} to this case which easily illustrates the main ideas of our method. In Section \ref{n}, we treat the case of arbitrary space dimension where a pseudodifferential change of variable will be needed, leading to a more technical proof. 

\smallskip


\section{Construction of the regularised problem}\label{rp}

For any given Schwartz function $\phi \in \mathscr{S}(\R^{n})$ we recall the notation 
$$
\phi_{\varepsilon} (x) := \frac{1}{\varepsilon^{n}} \phi\left(\frac{x}{\varepsilon}\right), \quad \varepsilon \in (0,1].
$$ 

\subsection{Regularisation of $a(t)$}

Let $\phi \in \mathscr{S}(\R)$ satisfying $0 \leq \phi \leq 1$ and $\int \phi = 1$. In order to regularise the coefficient $a(t)$ we first extend it to the whole real line, then we convolve it with the family $\phi_{\varepsilon}$ and lastly we restrict the result to the interval $[0,T]$. More to the point, we consider the extension 
$$
\tilde{a}(t) = 
\begin{cases}
	a(t), \quad\,\, t \in [0,T], \\
	a(T), \quad t \in [T,T+1], \\
	a(0), \quad\, t \in [-1,0], \\
	0, \quad\quad\,\,\, t \in \R - [-1,T+1].	
\end{cases}
$$
Then we set $\tilde{a}_{\varepsilon} := \tilde{a} \ast \phi_{\varepsilon}$, $\varepsilon \in (0,1]$, and 
\begin{equation}\label{eq_regularization_a}
a_{\varepsilon} := \left.\tilde{a}_{\varepsilon}\right|_{[0,T]}.
\end{equation}

Now we discuss a uniform lower and upper bound for $a_{\varepsilon}$ with respect to the parameter $\varepsilon$. By definition, for all $t \in [0,T]$, we have 
	$$
	a_{\varepsilon}(t) = \tilde{a}_{\varepsilon}(t) = \int_{\R} \tilde{a}(s) \phi_{\varepsilon}(t-s) ds = \int_{-1}^{T+1}  a(s) \phi_{\varepsilon}(t-s) ds.
	$$
Since $C_{a} \leq a(s) \leq \tilde{C}_{a}$ we get 
	$$
	\frac{C_{a}}{\varepsilon} \int_{-1}^{T+1} \phi \left(\frac{t-s}{\varepsilon}\right) ds  \leq a_{\varepsilon}(t) 
	\leq \frac{\tilde{C}_{a}}{\varepsilon} \int_{-1}^{T+1} \phi \left(\frac{t-s}{\varepsilon}\right) ds.
	$$
A suitable change of variables implies
	$$
	C_{a} \int_{-\frac{T+1-t}{\varepsilon}}^{\frac{t+1}{\varepsilon}} \phi(s) ds  \leq a_{\varepsilon}(t) 
	\leq \tilde{C}_{a} \int_{-\frac{T+1-t}{\varepsilon}}^{\frac{t+1}{\varepsilon}} \phi(s) ds.
	$$
Hence, taking $\varepsilon < \varepsilon_{0}$, where $\varepsilon_{0}$ is small enough, we get
	\begin{equation}\label{eq_uniform_bounds_reg_a}
	\frac{C_{a}}{2} \leq a_{\varepsilon}(t) \leq \tilde{C}_{a}, \quad t \in [0,T],\, \varepsilon \in (0, \varepsilon_{0}).
	\end{equation}
Indeed, denoting $I_{\varepsilon}(t) = \int_{-\frac{T+1-t}{\varepsilon}}^{\frac{t+1}{\varepsilon}} \phi(s) ds$, it suffices to notice that $\lim_{\varepsilon \to 0^{+}} I_{\varepsilon}(t) = 1$ uniformly in $t \in [0,T]$.

\subsection{Regularisation of $a_j(t)$, $j=0, 1, \ldots, n$} Let $\varphi \in \mathscr{S}(\R)$ such that $\int \varphi = 1$. We then define
\begin{equation}\label{eq_regularization_a_j}
	a_{j,\varepsilon} (t) = (a_{j} \ast \varphi_{\varepsilon}) (t) = \frac{1}{\varepsilon} a_j\left( \varphi \left( \frac{t-\cdot}{\varepsilon} \right) \right).
\end{equation}
Since $a_j\in \mathcal{E}'(I)\subset\mathcal{E}'(\R)$, $j = 0, 1, \ldots, n$, there exist a compact set $K \subset \R$, $C > 0$ and $\tilde{N}_0 \in \N_0$ such that
$$
|a_{j}(h)| \leq C \sum_{\alpha \leq \tilde{N}_0} \sup_{x \in K} |\partial^{\alpha}_{x}h(x)|, \quad \forall \, h \in C^{\infty}(\R).
$$
Hence, we immediately get
\begin{equation}\label{eq_upper_bound_coefficients_a_j}
	|a_{j,\varepsilon} (t)| \leq C \varepsilon^{-N_0} \max_{|\alpha| \leq N_0} \| \partial^{\alpha} \varphi\|_{\infty},
\end{equation}
where $N_0$ is the maximum of the orders of $a_j$, $j =0, 1, \ldots, n$, plus $1$.

\subsection{Regularisation of $b_j(x)$, $j=0,1,\ldots,n$} 
Let $\rho \in \mathscr{S}(\R^n)$ with $\int \rho = 1$. We then define 
\begin{equation}\label{eq_estimates_regularized_b_1_0}
	b_{j,\varepsilon}(x) := (b_j\ast \rho_{\varepsilon}) =  \frac{1}{\varepsilon^{n}} b_j\left( \rho \left( \frac{x-\cdot}{\varepsilon} \right) \right).
\end{equation}
Since $b_j \in H^{-\infty,2}(\R^{n})$, $j =1, \ldots,n$, and $b_0 \in H^{-\infty,0}(\R^{n})$, we get $b_{j,\epsilon} \in \mathcal{B}^{\infty}(\R^{n})$, $j =0,1,\ldots,n$, and the following estimates hold
\begin{equation}\label{regbj}
	|\partial^{\beta}_{x}b_{j,\varepsilon}(x)| \leq C_{\beta} \varepsilon^{-|\beta|-N_1} \langle x \rangle^{-2}, \quad
	|\partial^{\beta}_{x}b_{0,\varepsilon}(x)| \leq C_{\beta} \varepsilon^{-|\beta|-N_1},
\end{equation}
where $N_1 > 0$ is a number depending on the coefficients $b_j$, $j = 0, 1, \ldots, n$, and on the dimension.

\subsection{Regularisation of the Cauchy data} 
Let $\mu \in \mathscr{S}(\R^{n})$ with $\int \mu = 1$.\\ Let  $f \in C([0,T], H^{-\infty, n+1}(\R^{n}))$ and  $g\in H^{-\infty, n+1}(\R^{n})$.
We then define $f_\epsilon(t,x) = \mu_{\varepsilon}(x) \ast f(t,x)$ and $g_\epsilon(x) = \mu_{\varepsilon}(x) \ast g(x)$. The following estimates hold:
\begin{equation}
	|\partial^{\beta}_{x}f_{\varepsilon}(t,x)|\leq C_{\beta} \varepsilon^{-|\beta|-\tilde{N}_f} \langle x \rangle^{-(n+1)},\quad 
	|\partial^{\beta}_{x} g_\epsilon| \leq C_{\beta} \varepsilon^{-|\beta|-\tilde{N}_g} \langle x \rangle^{-(n+1)},
\end{equation}
where $\tilde{N}_f > 0$ is a number depending on $f(t)$ and on the dimension $n$, $\tilde{N}_g > 0$ is a number depending on $g$ and on the dimension $n$. By these estimates we immediately get that for all $m \in \N_0$ there exist $C>0$ and $N_f\in\N_0$ and $N_g\in\N_0$ such that 
\[
\begin{split}
\Vert f_\varepsilon(t,\cdot)\Vert_{H^m}&\le C_{m,N_f}\varepsilon^{-N_f},\quad \forall t\in[0,T],\ \varepsilon\in (0,1]\\
\Vert g_\varepsilon\Vert_{H^m}&\le C_{m,N_g}\varepsilon^{-N_g},\quad \forall  \varepsilon\in (0,1].
\end{split}
\]
\medskip

We are finally ready to define the family of regularised Cauchy problems that we shall study in the next sections. We consider the family of regularised operators
\begin{equation}\label{regularized_operator}
	S_{\varepsilon} = D_t - a_{\varepsilon}(t)\triangle_x + \sum_{j=1}^{n}a_{j,\varepsilon}(t) b_{j,\varepsilon}(x) D_{x_j} + a_{0,\varepsilon}(t)b_{0,\varepsilon}(x), \quad t \in [0,T], x \in \R^{n}, \varepsilon \in (0,\varepsilon_0],
\end{equation}
and then the family of ragularised problems
\begin{equation}\label{regularized_CP_true}
	\begin{cases}
		S_{\varepsilon} v(t,x) = f(t,x), \quad t \in [0,T],\, x \in \R^{n}, \\
		v(0,x) = g(x), \quad \quad \,\, \, x \in \R^{n},
	\end{cases}
\end{equation}
where $f \in C([0,T];H^{\infty}(\R^{n}))$ and $g \in H^{\infty}(\R^{n})$. In the next sections we will obtain a net of solutions $(u_{\varepsilon})_{\varepsilon \in (0,\varepsilon_{0}]}$ where for every $\varepsilon$ the function $u_{\varepsilon} \in C([0,T]; H^{\infty}(\R^{n}))$  is the unique solution of the Cauchy Problem \eqref{regularized_CP_true}. Moreover, we will also derive energy estimates for the solutions $u_{\varepsilon}$ expliciting how the constants depend on the parameter $\varepsilon$. Finally, thanks to these energy estimates, we will be able to prove the main result of the paper.

\smallskip


\section{Pseudodifferential operators}\label{Mein_teil}

In this section we collect some results and definitions concerning pseudodifferential operators that we will employ in the next sections. For the proofs we address the reader to \cite{Kumano-Go}.

\begin{Def}
	Given $m\in\R$, $\rho \in [0,1], \delta \in [0,1)$, we denote by $S^{m}_{\rho,\delta}(\R^{n})$ the space of all smooth functions $p(x,\xi) \in C^{\infty}(\R^{2n})$ such that for any $\alpha, \beta \in \N^{n}_{0}$ there exists a positive constant $C_{\alpha,\beta}$ satisfying
	$$
	|\partial^{\alpha}_{\xi} \partial^{\beta}_{x} p(x,\xi)| \leq C_{\alpha,\beta} \langle \xi \rangle^{m+\delta|\beta|-\rho|\alpha|}.
	$$
	The Frech\'et topology of the space $S^{m}_{\rho,\delta}(\R^{2n})$ is induced by the following family of seminorms
	$$
	|p|^{(m)}_\ell := \max_{|\alpha| \leq \ell, |\beta| \leq \ell} \sup_{x, \xi \in \R^{2n}} |\partial^{\alpha}_{\xi} \partial^{\beta}_{x} p(x,\xi)| \langle \xi \rangle^{-m - \delta|\beta|+\rho|\alpha|}, \quad p \in S^{m}_{\rho,\delta}(\R^{2n}), \, \ell \in \N_0.
	$$
\end{Def}

\begin{Rem}
	When $\rho = 1$ and $\delta = 0$ we simply write $S^{m}(\R^{2n})$ isntead of $S^{m}_{1,0}(\R^{2n})$. When $\rho = \delta = 0$ we have $S^{0}_{0,0}(\R^{2n}) = \mathcal{B}^{\infty}(\R^{2n})$.
\end{Rem}

Given a symbol $p(x,\xi)$  we associate to it with the continuous operator on the Schwartz space of rapidly decreasing functions $p(x,D): \mathscr{S}(\R^{n}) \to \mathscr{S}(\R^{n})$ defined by 
$$
p(x,D) u(x) = \int e^{i\xi x} p(x,\xi) \widehat{u}(\xi) \dslash\xi, \quad u \in \mathscr{S}(\R^{n}).
$$
We will sometimes also use the notation $p(x,D) = op(p(x,\xi))$. The next result concerns the action of such operators in Sobolev spaces.

\begin{Th}\label{theorem_Calderon_Vaillancourt}[Calder\'on-Vaillancourt]
	Let $p \in S^{m}_{\rho,\delta}(\R^{2n})$. Then for any real number $s \in \R$ there exist $\ell := \ell(s,m) \in \N_0$ and $C:= C_{s,m} > 0$ such that 
	\begin{equation*}
		\| p(x,D)u \|_{H^{s}} \leq C |p|^{(m)}_{\ell} \| u \|_{H^{s+m}}, \quad \forall \, u \in H^{s+m}.
	\end{equation*}
	Moreover, when $m = s = 0$ we can replace $|p|^{(m)}_{\ell}$ by
	\begin{equation*}
		\max_{|\alpha| \leq \ell_1, |\beta| \leq \ell_2} \sup_{x, \xi \in \R^{n}} |\partial^{\alpha}_{\xi} \partial^{\beta}_{x} p(x,\xi)| \langle \xi \rangle^{-\delta|\beta| +\rho|\alpha| },
	\end{equation*}
	where 
	$$
	\ell_1 = 2\left\lfloor\frac{n}{2} +1\right\rfloor, \quad \ell_2 = 2\left\lfloor \frac{n}{2(1-\delta)} + 1 \right\rfloor.
	$$
\end{Th}

Now we consider the algebra properties of $S^{m}_{\rho,\delta}(\R^{n})$ with respect to the composition of operators. In the sequel $Os-$ in front of the integral sign stands for oscillatory integral. Let $p_j \in S^{m_j}_{\rho,\delta}(\R^{n})$, $j = 1, 2$, and define 
\begin{align}\label{eq_symbol_of_composition}
	q(x,\xi) &= Os- \iint e^{-iy \eta} p_1(x,\xi+\eta)p_2(x+y,\xi) dy \dslash\eta \\
	&= \lim_{\mu \to 0} \iint e^{-iy\eta} p_1(x,\xi+\eta)p_2(x+y,\xi) e^{-\mu^2|y|^2} e^{-\mu^2|\eta|^2} dy \dslash\eta. \nonumber
\end{align} 

\begin{Th}\label{Hotel_california}
	Let $p_j \in S^{m_j}_{\rho,\delta}(\R^{2n})$, $j = 1, 2$, and consider $q$ defined by \eqref{eq_symbol_of_composition}. Then $q \in S^{m_1+m_2}_{\rho,\delta}(\R^{2n})$ and $q(x,D) = p_1(x,D) p_2(x,D)$. Moreover, the symbol $q$ has the following asymptotic expansion
	\begin{align*}
		q(x,\xi) = \sum_{|\alpha| < N} \frac{1}{\alpha!} \partial^{\alpha}_{\xi}p_{1}(x,\xi)D^{\alpha}_{x}p_2(x,\xi) + r_N(x,\xi), 
	\end{align*} 
	where 
	$$
	r_N(x,\xi) = \sum_{|\gamma| = N} \frac{N}{\gamma!} \int_{0}^{1} (1-\theta)^{N-1} \, Os - \iint e^{-iy\eta} \partial^{\gamma}_{\xi} p_1(x,\xi+\theta\eta) D^{\gamma}_{x} p_2(x+y,\xi)  dy\dslash\eta \, d\theta,
	$$
	and the seminorms of $r_N$ may be estimated in the following way: for any $\ell_{0} \in \N_0$ there exists $C_{\ell, N, n} > 0$ such that 
	$$
	|r_N|^{(m_1+m_2)}_{\ell_0} \leq C_{\ell,N,n} |p_1|^{(m_1)}_{\ell_0+N+n+1} |p_2|^{(m_2)}_{\ell_0+N+n+1}.
	$$
\end{Th}

The last theorem that we recall is the celebrated sharp G{\aa}rding inequality (see Theorem 2.1.3 in \cite{Kenig_SG}).

\begin{Th}\label{Theorem_sharp_garding}
	Let $p \in S^{1}(\R^{2n})$ and suppose $Re\,p(x,\xi) \geq 0$ for all $x \in \R^n$ and $|\xi| \geq R$ for some $R>0$. Then there exist $k = k(n) \in \N_0$ and $C = C(n,R)$ such that 
	$$
	Re\, \langle p(x,D) u, u \rangle_{L^{2}} \geq -C |p|^{(1)}_{k} \|u\|^{2}_{L^{2}}, \quad u \in \mathscr{S}(\R^{n}).
	$$ 
\end{Th}

\smallskip


\section{Solving the regularised problem: the case $n=1$}\label{n=1}
Let us consider the Cauchy problem \eqref{CP} in the case $n=1$. We apply to the coefficients and to the Cauchy data the regularisation described in Section \ref{rp} and we come to the regularised problem \eqref{regularized_CP_true} for the operator
\begin{equation}\label{regularized_operator1}
	S_{\varepsilon} = D_t + a_{\varepsilon}(t)D_x^2 + a_{1,\varepsilon}(t) b_{1,\varepsilon}(x) D_{x} + a_{0,\varepsilon}(t)b_{0,\varepsilon}(x), \quad t \in [0,T], x \in \R, \varepsilon \in (0,\varepsilon_0].
\end{equation}
We now consider the function
$$
B_{1,\varepsilon} (x) = \int_{0}^{x} b_{1,\varepsilon}(y) dy, \quad x \in \R.
$$
Using the fact  that $|b_{1,\varepsilon}| \leq C \varepsilon^{-N_1} \langle x \rangle^{-2}$ we conclude $B_{1,\varepsilon} \in \mathcal{B}^{\infty}(\R)$ and 
\begin{equation}\label{Bbounded}
|\partial^{\beta}_{x}B_{1,\varepsilon}(x)| \leq C_{\beta} \varepsilon^{-N_1 - \beta}.
\end{equation}
Since $a_{\varepsilon}$ never vanishes, we may define 
$$
e^{F_{\varepsilon}(t,x)}, \quad F_{\varepsilon}(t,x) = \frac{ia_{1,\varepsilon}(t) B_{1,\varepsilon}(x)}{2a_{\varepsilon}(t)} , \quad t \in [0,T],\, x \in \R.
$$
Conditions \eqref{eq_uniform_bounds_reg_a} and \eqref{Bbounded} imply that $F_{\varepsilon} \in C([0,T];\mathcal{B}^{\infty}(\R))$ and
\begin{equation}
	\left|e^{F_{\varepsilon}(t,x)}\right| \leq \text{exp}\left(C C_{a}^{-1}\varepsilon^{-N_0 - N_1}\right).
\end{equation}

Our idea, inspired by \cite{mizohata2014}, is to eliminate the first order coefficient by performing the following conjugation 
$$
S_{\varepsilon, F_{\varepsilon}}:= e^{F_{\varepsilon}(t,x)} \circ S_{\varepsilon} \circ e^{-F_{\varepsilon}(t,x)}.
$$
We have: 
\begin{align*}
	e^{F_{\varepsilon}} \circ D_{t} \circ e^{-F_{\varepsilon}} = D_t - D_{t}F_{\varepsilon} = D_t  + \frac{B_{1,\varepsilon}(x)}{2} \cdot \frac{a_{1,\varepsilon}(t)a'_{\varepsilon}(t) - a'_{1,\varepsilon}(t)a_{\varepsilon}(t)}{a_{\varepsilon}^{2}(t)},
\end{align*}
\begin{align*}
	e^{F_{\varepsilon}} \circ a_{\varepsilon}(t) D^{2}_{x} \circ e^{-F_{\varepsilon}} &= a_{\varepsilon}(t) D^{2}_{x} - a_{1,\varepsilon}(t)b_{1,\varepsilon}(x)D_{x} + a_{\varepsilon}(t)\{\partial^{2}_{x}F_{\varepsilon} - (\partial_{x}F_{\varepsilon})^{2}\} \\
	&=  a_{\varepsilon}(t) D^{2}_{x} - a_{1,\varepsilon}(t)b_{1,\varepsilon}(x)D_{x} + \frac{i}{2} a_{1,\varepsilon}(t) \partial_{x} b_{1,\varepsilon}(x) + \frac{1}{4a_{\varepsilon}(t)} (a_{1,\varepsilon}(t)b_{1,\varepsilon}(x))^{2},
\end{align*}
\begin{align*}
	e^{F_{\varepsilon}} \circ a_{1,\varepsilon}(t)b_{1,\varepsilon}(x) D_{x} \circ e^{-F_{\varepsilon}} &= a_{1,\varepsilon}(t) b_{1,\varepsilon}(x) D_{x} - \frac{(a_{1,\varepsilon}(t)b_{1,\varepsilon}(x))^{2}}{2a_{\varepsilon}(t)}.
\end{align*}
Therefore 
\begin{align*}
	S_{\varepsilon, F_{\varepsilon}} = D_t + a_{\varepsilon}(t) D^{2}_{x} + A_{0,\varepsilon}(t,x),
\end{align*}
where 
\begin{align*}
	A_{0,\varepsilon}(t,x) &= a_{0,\varepsilon}(t)b_{0,\varepsilon}(x) +  \frac{B_{1,\varepsilon}(x)(a_{1,\varepsilon}(t)a'_{\varepsilon}(t) - a'_{1,\varepsilon}(t)a_{\varepsilon}(t))}{2a_{\varepsilon}^{2}(t)} 
	+ \frac{i}{2} a_{1,\varepsilon}(t) \partial_{x} b_{1,\varepsilon}(x) - \frac{(a_{1,\varepsilon}(t)b_{1,\varepsilon}(x))^{2}}{4a_{\varepsilon}(t)}.
\end{align*}
We remark that the operator $S_{\varepsilon, F_{\varepsilon}}$ maintains the same structure as \eqref{regularized_operator1} but has no terms of order one. 
Writing $\partial_t = iS_{\varepsilon, F_{\varepsilon}} - ia_{\varepsilon}(t) D^{2}_{x} - iA_{0,\varepsilon}(t,x)$, we obtain

\begin{align*}
	\partial_{t} &\|u(t)\|^{2}_{L^2} = 2 Re\, (\partial_t u(t), u(t))_{L^2} \\
	&= 2 Re\, (iS_{\varepsilon, F_\varepsilon} u(t), u(t))_{L^2} - 2 \underbrace{Re\, (ia_{\varepsilon}(t)D^{2}_{x} u(t), u(t))_{L^2}}_{=0} - 2 Re\, (iA_{0,\varepsilon}(t) u(t), u(t))_{L^2} \\
	&\leq \|S_{\varepsilon,F_{\varepsilon}}u(t)\|^{2}_{L^2} + \left(1 + 2\sup_{\overset{t \in [0,T]}{x \in \R} } |A_{0,\varepsilon}(t,x)|\right) \|u(t)\|^{2}_{L^2}.
\end{align*}
Applying Gronwall lemma we get
\begin{align*}
	\|u(t)\|^{2}_{L^2} &\leq \text{exp}\left(t\left(1+2\sup_{t,x}|A_{0,\varepsilon}(t,x)|\right)\right) \left\{ \|u(0)\|^{2}_{L^2} + \int_{0}^{t} \|S_{\varepsilon,F_{\varepsilon}}u(\tau)\|^{2}_{L^2} d\tau \right\}.
\end{align*}
We estimate $A_{0,\varepsilon}(t,x)$ in the following way 
$$
|A_{0,\varepsilon}(t,x)| \leq C C_{a}^{-2} \varepsilon^{-2(N_0+N_1)-1}
$$
where $C > 0$ denotes a constant which may depend on $a, a_1, a_0, b_1, b_0, \varphi, \phi, \rho$, but does not depend on $\varepsilon$.  Therefore 
\begin{align}\label{powerslave}
	\|u(t)\|^{2}_{L^2} &\leq \text{exp}\left(tC C^{-2}_{a} \varepsilon^{-2(N_0+N_1)-1}\right) \left\{ \|u(0)\|^{2}_{L^2} + \int_{0}^{t} \|S_{\varepsilon,F_{\varepsilon}}u(\tau)\|^{2}_{L^2} d\tau \right\}.
\end{align}

To obtain \eqref{powerslave} with general Sobolev norms, $m \in \N_0$, we observe 
$$
S_{\varepsilon, F_{\varepsilon}, m} := \langle D_x \rangle^{m} \circ S_{\varepsilon, F_{\varepsilon}} \circ \langle D_x \rangle^{-m} = D_t + a_{\varepsilon}(t)D^{2}_{x} + \langle D_x \rangle^{m}A_{0,\varepsilon}(t,x) \langle D_x \rangle^{-m}.
$$
We have 
$$
\langle D_x \rangle^{m}A_{0,\varepsilon}(t,x) \langle D_x \rangle^{-m} = \sigma_{m,\varepsilon} (t,x,D)
$$
where 
\begin{align*}
	\sigma_{m,\varepsilon}(t,x,\xi) &= Os- \iint e^{-iy\eta} \langle \xi +\eta \rangle^{m} A_{0,\varepsilon}(t,x+y) \langle \xi \rangle^{-m} dy\dslash\eta \\
	&= \iint e^{-iy\eta} \langle y \rangle^{-2} \langle D_\eta \rangle^{2} \left\{ \langle \eta \rangle^{-2\lceil \frac{m+1}{2} \rceil} \langle D_y \rangle^{2\lceil \frac{m+1}{2} \rceil} \left( \langle \xi +\eta \rangle^{m} A_{0,\varepsilon}(t,x+y) \langle \xi \rangle^{-m} \right)   \right\}.
\end{align*}
Calder\'on-Vaillancourt Theorem implies  
$$
\|\sigma_{m,\varepsilon}(t,x,D)u\|_{L^2} \leq C_{m} C^{-2}_{a} \varepsilon^{-2(N_1+N_0) - m - 5} \| u \|_{L^2}. 
$$

In this way, proceeding in analogous manner as before, we obtain for any $m \in \N_0$ the following a-priori estimate:
\begin{align}\label{aces_high1}
	\|u(t)\|^{2}_{H^{m}} &\leq \text{exp}\left(tC_m C^{-2}_{a} \varepsilon^{-2(N_0+N_1)-m-5}\right) \left\{ \|u(0)\|^{2}_{H^{m}} + \int_{0}^{t} \|S_{\varepsilon,F_{\varepsilon}}u(\tau)\|^{2}_{H^{m}} d\tau \right\}.
\end{align}

Energy inequality \eqref{aces_high1} gives the following proposition.

\begin{Prop}
	Let $\tilde{f} \in C([0,T]; H^m(\R))$ and $\tilde{g} \in H^{m}(\R)$. There exists a unique solution $u$ in $C([0,T];H^{m}(\R))$ for the Cauchy problem 
	\begin{equation}\label{sign_of_the_cross}
		\begin{cases}
			S_{\varepsilon,F_{\varepsilon}} u(t,x) = \tilde{f}(t,x), \quad t \in [0,T], x \in \R, \\
			u(0,x) = \tilde{g}(x), \qquad \quad \,\,\, 	x \in \R,
		\end{cases}
	\end{equation}
	and the solution $u$ satisfies \eqref{aces_high1}.
\end{Prop}

The above proposition gives $H^{\infty}$-well-posedness for the problem \eqref{regularized_CP_true}. Indeed, first let us consider the auxiliary problem 
\begin{equation}\label{finocchio1}
	\begin{cases}
		S_{\varepsilon, F_{\varepsilon}} v(t,x) = e^{F_{\varepsilon}(t,x)}f(t,x), \quad t\in [0,T], x \in \R, \\
		v(0,x) = e^{F_{\varepsilon}(0,x)}g(x), \quad \qquad \,\,\, x \in \R.
	\end{cases}
\end{equation}
It is not hard to conclude that for $h\in H^m(\R)$ ($m \in \N_0$)
\begin{align*} 
\|e^{\pm F_{\varepsilon}(t,\cdot)} h(\cdot)\|_{H^m} &\leq C_{m} \max_{\beta \leq m} \sup_{x\in\R} |\partial^{\beta}_{x}e^{\pm F_{\varepsilon}(t,x)}| \,\|h\|_{H^{m}} \\
&\leq C_{m} \varepsilon^{-m} e^{C\varepsilon^{-(N_0+N_1)}} \|h\|_{H^{m}}.
\end{align*}

Let $v_{\varepsilon} \in C([0,T], H^{\infty}(\R))$ be the solution of \eqref{finocchio1}. Then $u_{\varepsilon} = e^{-F_{\varepsilon}}v_{\varepsilon}$ defines a solution for the regularised Cauchy problem \eqref{regularized_CP_true}. Next we observe that the energy inequality \eqref{aces_high1} implies for every $m\in \N_0$
\begin{align*}
	\|&u_{\varepsilon}(t)\|^{2}_{H^m} = \| e^{-F_{\varepsilon}(t)} v_{\varepsilon}(t)\|^{2}_{H^m} \leq C_m \varepsilon^{-2m}e^{2C\varepsilon^{-2(N_0+N_1)}} \|v_{\varepsilon}(t)\|^{2} \\
	&\leq C_m \varepsilon^{-2m}e^{2C\varepsilon^{-2(N_0+N_1)}} \text{exp}\left(tC_m \varepsilon^{-2(N_0+N_1)-m-5}\right) \left\{ \|e^{F_{\varepsilon}(0)}g\|^{2}_{H^{m}} + \int_{0}^{t} \|e^{F_{\varepsilon}(\tau)}f(\tau)\|^{2}_{H^{m}} d\tau \right\} \\
	&\leq C_m \varepsilon^{-4m}e^{4C\varepsilon^{-2(N_0+N_1)}} \text{exp}\left(tC_m \varepsilon^{-2(N_0+N_1)-m-5}\right) \left\{ \|g\|^{2}_{H^{m}} + \int_{0}^{t} \|f(\tau)\|^{2}_{H^{m}} d\tau \right\}.
\end{align*}
The uniqueness of the found solution follows by standard arguments. We summarise what we have done in the following theorem.

\begin{Th}\label{th}\label{appr1}
	For every $\varepsilon \in (0,\varepsilon_{0}]$, where $\varepsilon_0$ is a small parameter so that \eqref{eq_uniform_bounds_reg_a} holds, consider the regularised Cauchy problem \eqref{regularized_CP_true} with initial data $f \in C([0,T];H^{\infty}(\R))$ and $g \in H^{\infty}(\R)$. Then there exists a unique solution $u_{\varepsilon} \in C([0,T];H^{\infty}(\R))$ for the problem \eqref{regularized_CP_true}. Besides, the solution $u_\varepsilon$ satisfies for every $m\in\N_0$ 
	\begin{equation}\label{national_acrobat1}
		\|u_{\varepsilon}\|^{2}_{H^m} \leq C_{m}\text{exp}\left(C_{m,T} \varepsilon^{-2(N_0+N_1)-\theta_{m}}\right) \left\{ \|g\|^{2}_{H^{m}} + \int_{0}^{t} \|f(\tau)\|^{2}_{H^{m}} d\tau \right\},
	\end{equation}
	where 
	\begin{itemize}
		\item[(i)] $C_{m}$ and $C_{m,T}$ are constants depending on the coefficients $a,a_1,a_0,b_1,b_0$ and on the mollifiers $\rho, \phi, \varphi$; 
		
		\item [(ii)] $N_0$ stands for the maximum of the orders of $a_0$ and $a_1$ plus $1$;
		
		\item[(iii)] $N_1$ is a positive number depending on the coefficients $b_0, b_1$;
		
		\item[(iv)] $\theta_{m}$ is a natural number depending on $m$.
	\end{itemize}
\end{Th}

\smallskip


\section{Solving the regularised problem: the $n$-dimensional case}\label{n}

Let us now come to the regularised problem \eqref{regularized_CP_true} for the operator \eqref{regularized_operator} in the case of space dimension $n\geq 2$. To solve this problem we need  a change of variable as in the monodimensional case but in this more general situation, the change cannot be driven by a function, it has to be expressed in terms of a pseudodifferential operator. 

We will follow the argument used in \cite{KB}.
 For any $M > 0$ there exists a real-valued function $\tilde\lambda (x,\xi)$ satisfying the following conditions
$$
|\partial^{\alpha}_{\xi} \partial^{\beta}_{x} \tilde\lambda(x,\xi)| \leq M C_{\alpha,\beta} |\xi|^{-|\alpha|}, \quad x \in \R^n, |\xi| \geq 1, \alpha,\beta \in \N^{n}_{0},
$$
$$
\sum_{j=1}^{n} \xi_j \partial_{x_j} \tilde\lambda(x,\xi) \leq -M \langle x \rangle^{-2} \chi\left( \frac{\langle x \rangle}{|\xi|} \right) |\xi|,
$$
where $\chi(t) \in C^{\infty}_{c}(\R)$ satisfies $\chi(t) = 1$ for $|t| \leq \frac{1}{2}$, $\chi(t) = 0$ for $|t| \geq 1$, $t \chi'(t) \leq 0$ and $0 \leq \chi (t) \leq 1$. For a large parameter $h \geq 1$ to be chosen later on, we consider 
\begin{equation}\label{lambdan}
\lambda(x,\xi) = \tilde\lambda(x,\xi) (1-\chi)\left( \frac{ |\xi| }{ h } \right).
\end{equation}
Then, since $\langle \xi \rangle_{h} \leq \sqrt{5}h$ on the support of $(1-\chi)(h^{-1}|\xi|)$,
$$
|\partial^{\alpha}_{\xi} \partial^{\beta}_{x} \lambda(x,\xi)| \leq M C_{\alpha,\beta} \langle \xi \rangle^{-|\alpha|}_{h},
$$
where $\langle \xi \rangle_{h} := \sqrt{h^2+|\xi|^{2}}$ and the constants $C_{\alpha,\beta}$ do not depend on $M$ and on $h$. So, $\lambda$ is a symbol of order zero with respect to the basic weight $\langle \xi \rangle_{h}$. Moreover, we also have 
\begin{equation}\label{estimatelambda}
\sum_{j=1}^{n} \xi_j \partial_{x_j} \lambda(x,\xi) \leq -M \langle x \rangle^{-2} \chi\left( \frac{\langle x \rangle}{|\xi|} \right) |\xi|(1-\chi)\left( \frac{ |\xi| }{ h } \right).
\end{equation}

For a suitable choice of $M$ and $h := h(M)$ ($M$ and $h$ will depend on the parameter $\varepsilon$) we will prove that $e^{\lambda}(x,D)$ is invertible and the conjugated operator 
$$
S_{\varepsilon, \lambda} := e^{\lambda}(x,D) \circ S_{\varepsilon} \circ \{ e^{\lambda}(x,D)\}^{-1}
$$
satisfies a priori energy estimates, yielding a well-posed associated Cauchy problem. The main effort here is to check carefully how the constants in the energy estimate depend on the parameters $M, h$ and $\varepsilon$.

\subsection{Invertibility of $e^{\lambda(x,D)}$}

We begin with the following elementary lemma. 

\begin{Lemma}\label{ziggy_stardust}
	The following estimates hold
	$$
	|e^{\pm \lambda(x,\xi)}| \leq e^{CM}, \quad |\partial^{\alpha}_{\xi}\partial^{\beta}_{x} e^{\pm \lambda(x,\xi)}| \leq M^{|\alpha+\beta|} C_{\alpha,\beta} \langle \xi \rangle^{-|\alpha|}_{h} e^{\pm \lambda(x,\xi)},
	$$
	where the constants $C_{\alpha,\beta}$ and $C$ do not depend on $M$ and on $h$. In particular, $e^{\pm\lambda}$ is a symbol of order zero and its seminorms can be estimated as 
	$$
	|e^{\pm\lambda}|^{(0)}_{\ell} \leq C_{\ell} M^{2\ell} e^{CM}, \quad \ell \in \N_0,
	$$ 
	where $C_{\ell}$ does not depend on $M$ and on $h$.
\end{Lemma}

Now we consider the composition 
$$
e^{\lambda}(x,D) \circ e^{-\lambda}(x,D) = \sigma(x,D),
$$
where 
\begin{align*}
	\sigma(x,\xi) &= Os - \iint e^{-iy\eta} e^{\lambda(x,\xi+\eta)} e^{-\lambda(x+y,\xi)} dy \dslash\eta \\
	&= 1 + \underbrace{ i\sum_{j=1}^{n} \partial_{\xi_j} \lambda(x,\xi)\partial_{x_j}\lambda(x,\xi) + r_{-2}(x,\xi)}_{:=r(x,\xi)},
\end{align*}
\begin{equation}\label{the_writing_on_the_wall}
r_{-2}(x,\xi) = \sum_{|\gamma|=2} \frac{2}{\gamma!} Os - \iint e^{-iy\eta} \int_{0}^{1} (1-\theta) \partial^{\gamma}_{\xi} e^{\lambda(x,\xi+\theta\eta)} d\theta D^{\gamma}_{x} e^{-\lambda(x+y,\xi)} dy\dslash\eta.
\end{equation}
The estimates provided by Lemma \ref{ziggy_stardust} give
$$
|r_{-2}|^{(-2)}_{\ell} \leq C_{\ell,n} M^{4\ell+2n+6}e^{2CM} \leq C_{\ell,n} e^{(2C+1)M}, \quad \ell \in \N_0,
$$
where $C_{\ell,n}$ and $C$ do not depend on $M$ and on $h$. Next we observe that
$$
e^{\lambda}(x,D) \circ e^{-\lambda}(x,D) = I + r(x,D),
$$
$r(x,\xi)$ has order $-1$ and 
\begin{align*}
|\partial^{\alpha}_{\xi}\partial^{\beta}_{x}r(x,\xi)| &\leq C_{\alpha,\beta,n} e^{(2C+1)M} \langle \xi \rangle^{-1 - |\alpha|}_{h} \\
&\leq C_{\alpha,\beta,n} e^{(2C+1)M} h^{-1} \langle \xi \rangle^{-|\alpha|}.
\end{align*}
Hence, if we chose $h \geq h_0 (M, n) := A e^{(2C+1)M}$, where $A > 0$ is a number depending on the dimension and on a finite number of derivates of $r(x,\xi)$, Calder\'on-Vaillancourt theorem implies
$$
r(x,D): L^2(\R^n) \to L^2(\R^n)
$$ continuously
and the norm $ \sup\{ \|r(x,D)u\|_{L^2} \colon \|u\| \leq 1 \} < 1$. Let us also notice that due to our choice of $h_0$, for all $h \geq h_0$ the zero seminorms of $r$ do not depend on $M$. So, 
$$
\{I+r(x,D)\}^{-1} = \sum_{j \geq 0} \{-r(x,D)\}^{j} : L^2(\R^n) \to L^2(\R^n)
$$
and Theorem I.1 of \cite{Kumano-Go} (pag. $372$) gives that 
$$
\{I+r(x,D)\}^{-1} = s(x,D)
$$
where $s(x,\xi)$ is a symbol of order zero satisfying: for any $\ell \in \N_0$ there exists $\ell' \in \N_0$ such that  
$$
|s|^{(0)}_{\ell} \leq C_{\ell} \{|r|^{(0)}_{\ell'}\}^{3\ell} \leq C_{\ell} C_{\ell',n}^{3\ell}.
$$

\begin{Lemma}\label{so_procuro_minha_paz}
	For all $h \geq h_0(M,n) := A e^{(2C+1)M}$ the operator $e^{\lambda}(x,D)$ is invertible and we have
	$$
	\left(e^{\lambda}(x,D)\right)^{-1} = e^{-\lambda}(x,D) \circ \sum_{j \geq 0} (-r(x,D))^{j}, 
	$$
	where, for $r_{-2}$ given by \eqref{the_writing_on_the_wall},
	$$
	r(x,\xi) = i\sum_{j=1}^{n} \partial_{\xi_j} \lambda(x,\xi)\partial_{x_j}\lambda(x,\xi) + r_{-2}(x,\xi).
	$$
	Moreover, 
	$$
	\sum_{j \geq 0} (-r(x,D))^{j} = s(x,D),
	$$
	where $s(x,\xi)$ is a zero order symbol and its zero seminorms do not depend on $h$ and on $M$.
\end{Lemma}

We close this subsection writing the symbol $s(x,\xi)$ in a convenient way:
\begin{align*}
	s(x, D) &= I - r(x,D) + \sum_{j \geq 2} (-r(x,D))^{j} = I - r(x,D) + (r(x,D))^{2} s(x,D)  \\
			&= I - op\left( i\sum_{j=1}^{n} \partial_{\xi_j} \lambda(x,\xi)\partial_{x_j}\lambda(x,\xi) \right) - r_{-2}(x,D) +  (r(x,D))^{2} s(x,D).
\end{align*}
So, 
\begin{equation}\label{hell_on_Earth}
	s(x,\xi) = 1 - i\sum_{j=1}^{n} \partial_{\xi_j} \lambda(x,\xi)\partial_{x_j}\lambda(x,\xi) + s_{-2}(x,\xi),
\end{equation}
where $s_{-2}(x,\xi)$ has order $-2$ and satifies 
$$
|\partial^{\alpha}_{\xi}\partial^{\beta}_{x} s_{-2}(x,\xi)| \leq C_{\alpha,\beta,n} e^{CM} \langle \xi \rangle^{-2-|\alpha|}_{h},
$$
where the constants $C_{\alpha,\beta,n}$ and $C$ do not depend on $M$ and on $h$.

\subsection{Computing the operator $S_{\varepsilon, \lambda}$}
   
We start computing the symbol of
$$
e^{\lambda}(x,D) D^{2}_{x_j} e^{-\lambda}(x,D) s(x,D),
$$
which is
\begin{align*}
 \xi^{2}_{j} + 2i\xi_j \partial_{x_j} \lambda(x,\xi) + i\sum_{\ell=1}^{n} \xi^{2}_{j} \partial_{\xi_\ell}\lambda(x,\xi)\partial_{x_\ell} \lambda(x,\xi) + q_{j,0}(x,\xi), 
\end{align*}
where
$$
q_{j,0}(x,\xi) = \sum_{|\gamma|=2} \frac{2}{\gamma!} Os-\iint e^{-iy\eta} \int_{0}^{1} (1-\theta) \{\partial^{\gamma}_{\xi} (e^{\lambda}\xi^{2}_{j})\}(x,\xi+\theta\eta) d\theta D^{\gamma}_{x}e^{-\lambda(x+y,\xi)} dy\dslash\eta,
$$
$$
|\partial^{\alpha}_{\xi}\partial^{\beta}_{x}q_{j,0}(x,\xi)| \leq C_{\alpha,\beta,n} e^{(2C+1)M} \langle \xi \rangle^{-|\alpha|}_{h}.
$$
Hence 
\begin{align}\label{down_payment_blues}
	e^{\lambda}(x,D) \circ a_{\varepsilon}(t) \sum_{j=1}^{n} D_{x_j}^{2} \circ & e^{-\lambda}(x,D) = \textrm{op} \left( a_{\varepsilon}(t) \sum_{j=1}^{n} \xi^{2}_{j}
	+ a_{\varepsilon}(t)2i\sum_{j=1}^{n}\xi_j \partial_{x_j} \lambda(x,\xi) \right) \\ \nonumber
	&+ \textrm{op} \left(a_{\varepsilon}(t)i \sum_{j=1}^{n} \xi^{2}_{j} \sum_{\ell=1}^{n} \partial_{\xi_\ell}\lambda(x,\xi)\partial_{x_\ell} \lambda(x,\xi) + a_{\varepsilon}(t) \sum_{j=1}^{n} q_{j,0}(x,\xi) \right).
\end{align}
Now note that 
\begin{align}\label{stairway_to_heaven}
	-i a_{\varepsilon}(t) \sum_{j=1}^{n} D_{x_j}^{2}\ \circ \sum_{\ell=1}^{n} (\partial_{\xi_\ell} \lambda \partial_{x_\ell}\lambda)(x,D) &= 
	\textrm{op} \left(-i a_{\varepsilon}(t) \sum_{j=1}^{n} \xi^{2}_{j} \sum_{\ell=1}^{n} \partial_{\xi_\ell} \lambda(x,\xi)\partial_{x_\ell}\lambda(x,\xi) \right) \\ \nonumber
	&-\textrm{op} \left(2a_{\varepsilon}(t) \sum_{j=1}^{n} \xi_{j} \sum_{\ell=1}^{n} \partial_{x_j}\{\partial_{\xi_\ell} \lambda(x,\xi)\partial_{x_\ell}\lambda(x,\xi)\} \right) \\ \nonumber
	&+ \textrm{op} \left(ia_{\varepsilon}(t) \sum_{j=1}^{n} \sum_{\ell=1}^{n} \partial^{2}_{x_j}\{\partial_{\xi_\ell} \lambda(x,\xi)\partial_{x_\ell}\lambda(x,\xi)\} \right).
\end{align}
From \eqref{down_payment_blues} and \eqref{stairway_to_heaven} we conclude that
\begin{align}\label{whole_lotta_rosie}
	e^{\lambda}(x,D) \circ a_{\varepsilon}(t) \sum_{j=1}^{n} D^{2}_{x_j} \circ \{e^{-\lambda}(x,D)\}^{-1} =
	a_{\varepsilon}(t) \sum_{j=1}^{n} D^{2}_{x_j} +  a_{\varepsilon}(t)2i\sum_{j=1}^{n} (\partial_{x_j}\lambda)(x,D) D_{x_j} + a_{\varepsilon}(t)q_{0}(x,D),
\end{align}
where 
\begin{align*}
	q_{0}(x,D) &=\textrm{op} \left( \sum_{j=1}^{n} q_{j,0}(x,\xi) -2 \sum_{j=1}^{n} \xi_{j} \sum_{\ell=1}^{n} \partial_{x_j}\{\partial_{\xi_\ell} \lambda(x,\xi)\partial_{x_\ell}\lambda(x,\xi)\} +  i \sum_{j=1}^{n} \sum_{\ell=1}^{n} \partial^{2}_{x_j}\{\partial_{\xi_\ell} \lambda(x,\xi)\partial_{x_\ell}\lambda(x,\xi)\} \right) \\ &+ 2i\textrm{op}\left(\sum_{j=1}^{n}\xi_j \partial_{x_j} \lambda(x,\xi) 
	+ i \sum_{j=1}^{n} \xi^{2}_{j} \sum_{\ell=1}^{n} \partial_{\xi_\ell}\lambda(x,\xi)\partial_{x_\ell} \lambda(x,\xi) + \sum_{j=1}^{n} q_{j,0}(x,\xi) \right) \\ &\circ \textrm{op}\left(-i\sum_{j=1}^{n} \partial_{\xi_j} \lambda(x,\xi)\partial_{x_j}\lambda(x,\xi)\right) \\
	&+ e^{\lambda}(x,D) \circ \sum_{j=1}^{n} D^{2}_{x_j} \circ e^{-\lambda}(x,D) \circ s_{-2}(x,D),
\end{align*}
\begin{equation}\label{whole_lotta_rosie_2}
|q_{0}(x,\xi)|^{(0)}_{\ell} \leq C_{\ell,n} e^{CM},
\end{equation}
for some positive constant $C$ independent from $\varepsilon, h$ and $M$.

\medskip

Next we study the conjugation of the first order terms. We have
\begin{align*}
	e^{\lambda}(x,D) \circ a_{j,\varepsilon}(t)b_{j,\varepsilon}(x) D_{x_j} \circ e^{-\lambda}(x,D) &= a_{j,\varepsilon}(t)e^{\lambda}(x,D) \circ \textrm{op}\left( b_{j,\varepsilon}(x) \xi_j e^{-\lambda(x,\xi)} +ib_{j,\varepsilon}(x)\partial_{x}\lambda(x,\xi)e^{-\lambda(x,\xi)}  \right) \\
	&= a_{j,\varepsilon}(t)b_{j,\varepsilon}(x) D_{x_j} + a_{j,\varepsilon}(t) e_{j,\varepsilon,0}(x,D),
\end{align*}
where
\begin{align*}
e_{j,\varepsilon,0} (x,\xi) &= Os- \iint e^{-iy\eta} e^{\lambda(x,\xi+\eta)} b_{j,\varepsilon}(x+y) \partial_x \lambda(x+y,\xi) e^{-\lambda(x+y,\xi)} dy\dslash\eta \\
&+  \sum_{|\gamma|=1} Os- \iint e^{-iy\eta} \int_{0}^{1} \partial^{\gamma}_{\xi} e^{\lambda(x,\xi+\theta\eta)}d\theta b_{j,\varepsilon}(x+y)\xi_j e^{-\lambda(x+y,\xi)} dy\dslash\eta.
\end{align*}
Setting 
\begin{align*}
e_{\varepsilon,0}(t,x,D) &= \sum_{j=1}^{n} a_{j,\varepsilon}(t) e_{j,\varepsilon,0}(x,D) \\ 
&- \sum_{j=1}^{n}\{ a_{j,\varepsilon}(t) b_{j,\varepsilon}(x) D_{x_j} - a_{j,\varepsilon}(t)e_{j,\varepsilon,0}(x,D)\} \circ \textrm{op}\left( i\sum_{j=1}^{n} \partial_{\xi_j} \lambda(x,\xi)\partial_{x_j}\lambda(x,\xi) - s_{-2}(x,\xi) \right),
\end{align*}
we get
\begin{equation}\label{Caught_somewhere_in_time}
	e^{\lambda(x,D)} \circ \sum_{j=1}^{n} a_{j,\varepsilon}(t)b_{j,\varepsilon}(x)D_{x_j} \circ \{e^{\lambda(x,D)}\}^{-1} = \sum_{j=1}^{n} a_{j,\varepsilon}(t) b_{j,\varepsilon}(x)D_{x_j} + e_{\varepsilon,0}(t,x,D),
\end{equation}
and we have the following estimate
\begin{equation}\label{Caught_somewhere_in_time_2}
	|\partial^{\alpha}_{\xi}\partial^{\beta}_{x} e_{\varepsilon,0}(t,x,\xi)| \leq C_{\alpha,\beta,n} \varepsilon^{-N_0-N_1-\ell_n-|\beta|} e^{CM}\langle \xi \rangle^{-|\alpha|}_{h},
\end{equation}
where $\ell_n$ is a natural number depending only on the dimension and $C_{\alpha,\beta, n}$, $C$ do not depend on $\varepsilon, M$ and $h$.

For the zero order term we have
\begin{align*}
	e^{\lambda}(x,D) \circ a_{0,\varepsilon}(t) b_{0,\varepsilon}(x) \circ \{e^{\lambda}(x,D)\}^{-1} = a_{0,\varepsilon}(t) \underbrace{\tilde{c}_{0,\varepsilon}(x,D) \circ s(x,D)}_{=:c_{0,\varepsilon}(x,D)}, 
\end{align*}
where 
$$
\tilde{c}_{0,\varepsilon}(x,\xi) = Os- \iint e^{-iy\eta} e^{\lambda(x,\xi+\eta)} b_{0,\varepsilon}(x+y) e^{-\lambda(x+y,\xi)} dy\dslash\eta.
$$
Hence 
\begin{equation}\label{hallowed_be_thy_name}
e^{\lambda} (x,D) \circ a_{0,\varepsilon}(t) b_{0,\varepsilon}(x) \circ \{e^{\lambda}(x,D)\}^{-1} = a_{0,\varepsilon}(t) c_{0,\varepsilon}(x,D),
\end{equation}
and we have the following estimate
\begin{equation}\label{hallowed_be_thy_name_2}
	|a_{0,\varepsilon}(t) \partial^{\alpha}_{\xi}\partial^{\beta}_{x} c_{0,\varepsilon}(x,\xi)| \leq C_{\alpha,\beta,n} \varepsilon^{-N_0-N_1-\iota_n-|\beta|} e^{CM}\langle \xi \rangle^{-|\alpha|}_{h},
\end{equation}
where $\iota_n$ is a natural number depending only on the dimension and $C_{\alpha,\beta, n}$, $C$ do not depend on $\varepsilon, M$ and $h$.

Combining \eqref{whole_lotta_rosie}, \eqref{whole_lotta_rosie_2}, \eqref{Caught_somewhere_in_time}, \eqref{Caught_somewhere_in_time_2}, \eqref{hallowed_be_thy_name} and \eqref{hallowed_be_thy_name_2} we conclude that
\begin{equation}\label{the_clasnman}
	S_{\varepsilon,\lambda} = D_t + a_{\varepsilon}(t)\sum_{j=1}^{n} D^{2}_{x_j} + \sum_{j=1}^{n} a_{j,\varepsilon}(t)b_{j,\varepsilon}(x)D_{x_j} + a_{\varepsilon}(t)2i\sum_{j=1}^{n} \partial_{x_j}\lambda(x,D) D_{x_j} + d_{0,\varepsilon}(t,x,D),
\end{equation}
where 
$$
d_{0,\varepsilon}(t,x,\xi) = a_{\varepsilon}(t)q_{0}(x,\xi) + e_{\varepsilon, 0} (t,x,\xi) + a_{0,\varepsilon}(t)c_{0,\varepsilon}(x,\xi),
$$
and we have the following estimate
\begin{equation}\label{the_clasnman_2}
	|\partial^{\alpha}_{\xi}\partial^{\beta}_{x} d_{0,\varepsilon}(t,x,\xi)| \leq C_{\alpha,\beta,n} \varepsilon^{-N_0-N_1-\theta_n-|\beta|} e^{CM}\langle \xi \rangle^{-|\alpha|}_{h},
\end{equation}
where $\theta_n$ denotes a natural number depending only on the dimension and the constants $C_{\alpha,\beta, n}$, $C$ do not depend on $\varepsilon, M$ and $h$.

\subsection{Energy estimate for the regularised problem.}

In order to derive an $L^2$ a priori energy estimate for $S_{\varepsilon,\lambda}$ we write 
\begin{align*}
iS_{\varepsilon,\lambda} &= \partial_t + ia_{\varepsilon}(t)\sum_{j=1}^{n} D^{2}_{x_j} + \sum_{j=1}^{n}\{-Im\,(a_{j,\varepsilon}(t) b_{j,\varepsilon}(x)) - 2a_{\varepsilon}(t)(\partial_{x_j}\lambda)(x,D) \}D_{x_j} \\
&+\frac{1}{2}\sum_{j=1}^{n} \{2i Re\,(a_{j,\varepsilon}(t)b_{j,\varepsilon}(x))D_{x_j} +\partial_{x_j}Re\,(a_{j,\varepsilon}(t)b_{j,\varepsilon}(x) )\} - \frac{1}{2} \sum_{j=1}^{n} \partial_{x_j} Re\,(a_{j,\varepsilon}(t)b_{j,\varepsilon}(x)) \\
 &+ d_{0,\varepsilon}(t,x,D).
\end{align*}

We immediately note that 
$$
Re\, \left\langle ia_{\varepsilon}(t)\sum_{j=1}^{n} D^{2}_{x_j} u, u \right\rangle_{L^2} = 0,
$$
$$
Re\, \left\langle \frac{1}{2}\sum_{j=1}^{n} \{2i Re\,(a_{j,\varepsilon}(t)b_{j,\varepsilon}(x))D_{x_j} +\partial_{x_j}Re\,(a_{j,\varepsilon}(t)b_{j,\varepsilon}(x) )\} u, u \right\rangle_{L^2} = 0
$$
and 
$$
 \left| \left\langle \frac{1}{2} \sum_{j=1}^{n} \partial_{x_j} Re\,(a_{j,\varepsilon}(t)b_{j,\varepsilon}(x)) u, u \right\rangle_{L^2} \right| \leq C_n \varepsilon^{-N_0 -N_1 -1} \|u\|^{2}_{L^2}.
$$
On the other hand, Calder\'on-Vaillancourt theorem implies 
$$
|\langle d_{0,\varepsilon}(t,x,D) u, u \rangle_{L^2}| \leq C_n \varepsilon^{-N_0 -N_1 - \theta_n} e^{CM} \|u\|^{2}_{L^2}.
$$

The estimate of 
$$
\left\langle  \sum_{j=1}^{n}\{-Im\,(a_{j,\varepsilon}(t) b_{j,\varepsilon}(x)) - 2a_{\varepsilon}(t)(\partial_{x_j}\lambda)(x,D) \}D_{x_j} u, u \right\rangle_{L^2}
$$
is more complicate. We first consider the split
\begin{align*}
-Im\,&(a_{j,\varepsilon}(t) b_{j,\varepsilon}(x)) \xi_j \{\chi(h^{-1}|\xi|) + (1-\chi)(h^{-1}|\xi|)\}  
\\
&= -Im\,(a_{j,\varepsilon}(t) b_{j,\varepsilon}(x)) \xi_j (1-\chi)(h^{-1}|\xi|)\chi(\langle x \rangle |\xi|^{-1}) 
\\
& \quad \quad  - Im\,(a_{j,\varepsilon}(t) b_{j,\varepsilon}(x)) \xi_j \chi(h^{-1}|\xi|) \chi(\langle x \rangle |\xi|^{-1}) 
 \\
 &\quad\quad \quad - Im\,(a_{j,\varepsilon}(t) b_{j,\varepsilon}(x)) \xi_j \chi(h^{-1}|\xi|) (1-\chi)(\langle x \rangle |\xi|^{-1}) \\
 &\quad\quad \quad\quad -Im\,(a_{j,\varepsilon}(t) b_{j,\varepsilon}(x)) \xi_j (1-\chi)(h^{-1}|\xi|) (1-\chi)(\langle x \rangle |\xi|^{-1}) .
\end{align*}
On the support of $(1-\chi)(\langle x \rangle |\xi|^{-1})$ it holds $\langle x \rangle^{-1} \leq |\xi|^{-1}$, so from the decay $\langle x \rangle^{-2}$ of the coefficients $b_{j,\varepsilon}$ we conclude
\begin{align*}
	|\partial^{\alpha}_{\xi} \partial^{\beta}_{x} \{Im\,(a_{j,\varepsilon}(t) b_{j,\varepsilon}(x)) \xi_j \chi(h^{-1}|\xi|) (1-\chi)(\langle x \rangle |\xi|^{-1})\}| \leq 
	C_{\alpha,\beta,n} \varepsilon^{-N_0-N_1-|\beta|} \langle \xi \rangle^{-|\alpha|}_{h},
\end{align*}
\begin{align*}
	|\partial^{\alpha}_{\xi} \partial^{\beta}_{x} \{Im\,(a_{j,\varepsilon}(t) b_{j,\varepsilon}(x)) \xi_j (1-\chi)(h^{-1}|\xi|) (1-\chi)(\langle x \rangle |\xi|^{-1})\}| \leq 
	C_{\alpha,\beta,n} \varepsilon^{-N_0-N_1-|\beta|} \langle \xi \rangle^{-|\alpha|}_{h}.
\end{align*}
Using the support properties of $\chi(h^{-1}|\xi|)$ we get $|\xi| \leq h$, hence
\begin{align*}
	|\partial^{\alpha}_{\xi} \partial^{\beta}_{x} \{Im\,(a_{j,\varepsilon}(t) b_{j,\varepsilon}(x)) \xi_j \chi(h^{-1}|\xi|) \chi(\langle x \rangle |\xi|^{-1})\}| \leq 
	C_{\alpha,\beta,n} h \varepsilon^{-N_0-N_1-|\beta|} \langle \xi \rangle^{-|\alpha|}_{h}.
\end{align*}

By the definition of the transformation $\lambda$ we get 
\begin{align}\nonumber
	p_{\varepsilon}(t,x,\xi) := \sum_{j=1}^{n}\{-&Im\,(a_{j,\varepsilon}(t) b_{j,\varepsilon}(x))(1-\chi)(h^{-1}|\xi|)\chi(\langle x \rangle |\xi|^{-1}) - 2a_{\varepsilon}(t)(\partial_{x_j}\lambda)(x,D) \}\xi_j \\\label{crucial point a(t)}
	&\geq  \{2C_{a}M - C(a_j,b_j,n)\varepsilon^{-N_0-N_1} \} |\xi| \langle x \rangle^{-2} (1-\chi)(h^{-1}|\xi|) \chi(\langle x \rangle |\xi|^{-1}).
\end{align}
Choosing 
$$
M = \frac{C(a_j,b_j,n)}{2C_a} \varepsilon^{-N_0-N_1}
$$
we obtain $p_{\varepsilon}(t,x,\xi) \geq 0$ for all $t,x$ and $\xi$. Besides, the following estimate holds
$$
|\partial^{\alpha}_{\xi}\partial^{\beta}_{x}p_{\varepsilon}(t,x,\xi)| \leq C_{n,\alpha,\beta} M \varepsilon^{-N_0-N_1-|\beta|} \langle \xi \rangle^{1-|\alpha|}_{h}. 
$$
Sharp G{\aa}rding inequality then gives 
$$
Re\, \langle p_{\varepsilon}(t,x,D) u, u \rangle \geq -C_{n} M\varepsilon^{-N_0-N_1-\theta_n} \|u\|^{2}_{L^{2}}.
$$

Finally, gathering all the computations above we get
\begin{align}\label{neon_knight}
	\partial_{t} \|u(t)\|^{2}_{L^2} &= 2 Re\, \langle \partial_t u, u \rangle_{L^2} \\
	&\leq \|S_{\varepsilon,\lambda} u\|^{2}_{L^2} + \|u\|^{2}_{L^2} + C_{n} \varepsilon^{-N_0-N_1-\theta_n} e^{CM} \|u\|^{2}_{L^2} \nonumber.
\end{align}
To obtain \eqref{neon_knight} with general Sobolev norms we apply the same argument to the operator 
$$
\langle D_x \rangle^{m} S_{\lambda,\varepsilon} \langle D_x \rangle^{-m}
$$
which is equal to $S_{\varepsilon,\lambda}$ plus a zero order term, which can be easily estimated using Theorem \ref{Hotel_california}. So, 
\begin{align*}
	\partial_{t} \|u(t)\|^{2}_{H^m} &\leq \|S_{\lambda,\varepsilon} u\|^{2}_{H^m} + \|u\|^{2}_{H^m} + C_{n,m} \varepsilon^{-N_0-N_1-\theta_{n,m}} e^{CM} \|u\|^{2}_{H^m},
\end{align*}
where $\theta_{n,m}$ denotes a natural number depending only on $n$ and $m$. Gronwall inequality and the definition of $M$ give
\begin{align}\label{aces_high}
	 \|u(t)\|^{2}_{H^m} &\leq \exp\left\{ C_{T,m,n} e^{\varepsilon^{-N_0-N_1-\theta_{n,m}}} \right\} \left( \|u(0)\|^{2}_{H^{m}} + \int_{0}^{t}\|S_{\varepsilon,\lambda}(\tau) u(\tau)\|^{2}_{H^m} d\tau \right).
\end{align}

	Let us summarise what we have done: 
	\begin{itemize}
		\item[-] first we proved that if
		$$
		h \geq h_0(M) = A e^{CM},
		$$ 
		for some constants $A, C > 0$ indepedent from $M$ and $\varepsilon$, then the operator $e^{\lambda}$ is invertible (c.f. Lemma \ref{so_procuro_minha_paz});
		
		\item[-] next, after some technical computations and estimates, we proved that if 
		$$
		M = \frac{C(a_j,b_j,n)}{2C_a} \varepsilon^{-N_0-N_1}
		$$
		then we may apply sharp G{\aa}rding inequality to get the energy inequality \eqref{aces_high}. 
		\end{itemize} In other words, setting $M$ and $h(M)$ as above we come to the energy inequality \eqref{aces_high} where the precise dependence on the parameter $\varepsilon$ is exhibited.

The next proposition is a consequence of \eqref{aces_high}.

\begin{Prop}
	Let $\tilde{f} \in C([0,T]; H^m(\R^{n}))$ and $\tilde{g} \in H^{m}(\R^{n})$. There exists a unique solution $u$ in $C([0,T];H^{m}(\R^{n}))$ to the Cauchy problem 
	\begin{equation}\label{sign_of_the_cross1}
		\begin{cases}
			S_{\varepsilon,\lambda} u(t,x) = \tilde{f}(t,x), \quad t \in [0,T], x \in \R^{n}, \\
			u(0,x) = \tilde{g}(x), \qquad \quad \, 	x \in \R^{n},
		\end{cases}
	\end{equation}
	and the solution $u$ satisfies \eqref{aces_high}.
\end{Prop}

We are finally ready to conclude the desired $H^{\infty}$ well-posedness for the problem \eqref{regularized_CP_true} in the $n$-dimensional case. We first consider the auxilliary problem 
\begin{equation}\label{finocchio}
	\begin{cases}
		S_{\varepsilon, \lambda} v(t,x) = e^{\lambda(x,D)}f(t,x), \quad t\in [0,T], x \in \R^{n}, \\
		v(0,x) = e^{\lambda(x,D)}g(x), \quad \qquad \, x \in \R^{n}.
	\end{cases}
\end{equation}
We observe that Calder\'on-Vaillancourt theorem gives
$$ 
	\|e^{\lambda}(x,D) u\|_{H^m} \leq C_{n,m} e^{CM} \|u\|_{H^{m}}, \qquad \|\{e^{\lambda}(x,D)\}^{-1} u\|_{H^m} \leq C_{n,m} e^{CM} \|u\|_{H^{m}},
$$
where $C_{n,m}$ and $C$ do not depend on $\varepsilon$, $M$ and $h$.

Let $v_{\varepsilon} \in C([0,T];H^{\infty}(\R^{n}))$ be the solution of \eqref{finocchio}. Then $u_{\varepsilon} = \{e^{\lambda(x,D)}\}^{-1}v_{\varepsilon}$ defines a solution for the regularised Cauchy problem \eqref{regularized_CP_true}. Next we observe that the energy inequality \eqref{aces_high} implies for every $m\in\N_0$
\begin{align*}
	\|&u_{\varepsilon}(t)\|^{2}_{H^m} = \| \{e^{\lambda}(\cdot,D)\}^{-1} v_{\varepsilon}(t)\|^{2}_{H^m} \leq C_{n,m} e^{CM} \|v_{\varepsilon}(t)\|^{2}_{H^{m}} \\
	&\leq C_{n,m} e^{CM} \exp\left\{ C_{T,m,n} e^{\varepsilon^{-N_0-N_1-\theta_{n,m}}} \right\} \left\{ \|e^{\lambda}(\cdot,D)g\|^{2}_{H^{m}} + \int_{0}^{t} \|e^{\lambda}(\cdot,D)f(\tau)\|^{2}_{H^{m}} d\tau \right\} \\
	&\leq C_{n,m}^{2} e^{2CM} \exp\left\{ C_{T,m,n} e^{\varepsilon^{-N_0-N_1-\theta_{n,m}}} \right\} \left\{ \|g\|^{2}_{H^{m}} + \int_{0}^{t} \|f(\tau)\|^{2}_{H^{m}} d\tau \right\}.
\end{align*}
The uniqueness of the solution follows by standard arguments. We summarise what we have done in the following theorem.

\begin{Th}\label{th}\label{apprn}
	For every $\varepsilon \in (0,\varepsilon_{0}]$, where $\varepsilon_0$ is a small parameter so that \eqref{eq_uniform_bounds_reg_a} holds, consider the regularised Cauchy Problem \eqref{regularized_CP_true} with initial data $f \in C([0,T];H^{\infty}(\R^{n}))$ and $g \in H^{\infty}(\R^{n})$. Then there exists a unique solution $u_{\varepsilon} \in C([0,T];H^{\infty}(\R^{n}))$ for the problem \eqref{regularized_CP_true}. Besides, the solution $u_\varepsilon$ satisfies for every $m\in\N_0$ the following estimate:
	\begin{equation}\label{national_acrobat}
		\|u_{\varepsilon}\|^{2}_{H^m} \leq  C_{n,m} \exp\left\{ C_{T,m,n} e^{\varepsilon^{-N_0-N_1-\theta_{n,m}}} \right\} \left\{ \|g\|^{2}_{H^{m}} + \int_{0}^{t} \|f(\tau)\|^{2}_{H^{m}} d\tau \right\},
	\end{equation}
	where 
	\begin{itemize}
		\item[(i)] $C_{n,m}$ and $C_{T,m,n}$ are constants depending on the coefficients $a,a_j,b_j$, $j=0,1,\ldots,n$, and on the mollifiers $\rho, \phi, \varphi$; 
		
		\item [(ii)] $N_0$ stands for the maximum of the orders of $a_0, a_1, \ldots, a_n$ plus $1$;
		
		\item[(iii)] $N_1$ is a positive number depending on the coefficients $b_0, b_1, \ldots, b_n$ and on the dimension;
		
		\item[(iv)] $\theta_{n,m}$ is a natural number depending only on $m$ and on the dimension $n$.
	\end{itemize}
\end{Th}
We remark that Theorem \ref{apprn} is the $n$-dimensional version of Theorem \ref{appr1} with a greater constant in the energy estimate. From now on we shall continue the proof of Theorem \ref{mainn} referring only to the $n$-dimensional case, since the final result for the regularised problem in the monodimensional case is a particular case.
\smallskip


\section{Proof of Theorem \ref{mainn}}\label{lady_evil}

We can now conclude the proof of Theorem \ref{mainn}. We know that the regularised Cauchy data $f_\epsilon (t), g_\epsilon$ fulfill the following estimate: for all $m\in\N_0$ there exists $C>0$, $N_f\in\N_0$ and $N_g\in\N_0$ such that 
\[
\begin{split}
\Vert f_\varepsilon(t,\cdot)\Vert_{H^m}&\le C\varepsilon^{-N_f},\\
\Vert g_\varepsilon\Vert_{H^m}&\le C\varepsilon^{-N_g},
\end{split}
\]
uniformly in $\varepsilon$ and $t$, and from Theorem \ref{apprn} the regularised problem \eqref{regularized_CP_true} with data $f_\varepsilon(t)$ and $g_\varepsilon$ has a unique solution $u_\varepsilon\in C([0,T];H^{\infty}(\R^{n}))$ satisfying for every $m\in\N_0$ 
\begin{align*}
\|u_{\varepsilon}\|^{2}_{H^m} &\leq  C_{n,m} \exp\left\{ C_{T,m,n} e^{\varepsilon^{-N_0-N_1-\theta_{n,m}}} \right\} \left\{ \|g_{\varepsilon}\|^{2}_{H^{m}} + \int_{0}^{t} \|f_{\varepsilon}(\tau)\|^{2}_{H^{m}} d\tau \right\} \\
	&\leq A_m  \,\text{exp}(B_me^{\varepsilon^{-N_m}}) \varepsilon^{-(N_f + N_g)}, \quad \varepsilon \in (0, \varepsilon_0), t \in [0, T]
\end{align*}
where $A_m, B_m > 0$ are independent from $\varepsilon$ and $N_m= N_0+N_1+\theta_{n,m} \in \N_0$ is independent from $\varepsilon$ but depending on the coefficients $b_j, a_j$, $j = 0, 1, \ldots, n$, on the dimension $n$ and on Sobolev index $m$.  

Next we can obtain moderate estimates for the net $(u_\varepsilon)_\varepsilon$ if we regularise the coefficients $b_j, a_j$, $j = 0,1, \ldots, n$, using the following positive scale
\begin{equation}\label{the_rise_of_evil}
\omega(\varepsilon) = \{ \log(\log(\log(\log(\varepsilon^{-1})))) \}^{-1}, \quad \varepsilon \in (0, \varepsilon_{0}).
\end{equation}
This means to replace $\varepsilon$ with $\omega(\varepsilon)$ in all the estimates above. First, we note that following elementary inequality: for any $r \geq 1$ we have
\begin{equation}\label{for_those_about_rock}
	\log(y) \leq C_{r} y^{\frac{1}{r}}, \quad \forall\, y \geq 2.
\end{equation}
Hence, denoting $X = \log(\log(\varepsilon^{-1}))$ we get
$$
\omega(\varepsilon)^{-N_m} = \{\log(\log(X)\}^{N_m} \leq C_{N_m}^{N_m} \log(X) = \log\left(X^{C^{N_m}_{N_m}}\right) \implies e^{\omega(\varepsilon)^{-N_m}} \leq X^{C^{N_m}_{N_m}}.
$$
Applying \eqref{for_those_about_rock} once more we conclude
$$
B_m X^{C^{N}_{N}} \leq B_m \tilde{C}_{N_m}^{C^{N_m}_{N_m}} \log(\varepsilon^{-1}) = \log(\varepsilon^{-D_{N_m}}),
$$
where $D_{N_m} \in \N_0$ is a large number depending on $N_m$ and on $B_m$. So, if use the regularisations
$$
 b_{j,\varepsilon} := \rho_{\omega(\varepsilon)} \ast b_j, \quad a_{j,\varepsilon} = \varphi_{\omega(\varepsilon)} \ast a_j, \quad j = 0, 1, \ldots, n,
$$
where $\omega$ is given by \eqref{the_rise_of_evil} we get the following estimate for the net of solutions $(u_\varepsilon)_{\varepsilon}$: for any $m \in \N_0$ we find constants $A_m, D_{N_m}$ such that 
$$
\|u_{\varepsilon}(t)\|^{2}_{H^m} \leq A_m \varepsilon^{-D_{N_m}-(N_f+N_g)}, \quad t \in [0,T], \, \varepsilon \in (0,\varepsilon_0). 
$$

By Definitions \ref{defmoderateness} and \ref{def_vw} we can conclude that the Cauchy problem \eqref{CP} admits very weak solutions of $H^\infty$-type. The main result of this paper is then proved.

\begin{Rem}
	We stress the fact that in order to obtain a $H^{\infty}$-moderate net of solutions we need to use the special positive scale $\omega$ (given in \eqref{the_rise_of_evil}) just for the coefficients of order one and order zero, meanwhile we are free to use the standard scale ($\omega(\varepsilon) = \varepsilon$) for the leading coefficient and Cauchy data.
\end{Rem}

\smallskip


\section{Uniqueness of the very weak solution and consistency with the regular theory}\label{the_tragic_prince}
We conclude the paper by discussing in which sense the found very weak solution for our Cauchy problem is unique and by proving that our result is consistent with the classical theory when the coefficients are regular enough.
\subsection{Uniqueness}
The very weak solution $(u_{\varepsilon})_{\varepsilon}$ that we found is unique in the following sense: if we perturb the coefficients of the regularised operator by suitable negligible nets, then the net of solution of the Cauchy problem associated with the perturbed operator will differ from $(u_{\varepsilon})_{\varepsilon}$ by a $H^{\infty}$-negligible net. Being more precise, suppose that we perturb the coefficients of the regularised operator $S_{\varepsilon}$ given by \eqref{regularized_operator1} in the following way
\begin{align*}
	S'_{\varepsilon} &:= D_t + (a_{\varepsilon}(t)+n_{\varepsilon}(t))\sum_{j=1}^{n} D^{2}_{x_j} \\
	&+ \sum_{j=1}^{n}(a_{j,\varepsilon}(t) +n_{j,1,\varepsilon}(t))(b_{j,\varepsilon}(x)+l_{j,1,\varepsilon}(x)) D_{x_j} + (a_{0,\varepsilon}(t)+n_{0,\varepsilon}(t))(b_{0,\varepsilon}(x)+l_{0,\varepsilon}(x))
\end{align*}
where 
\begin{itemize}
	\item[(a)] $n_{\varepsilon} \in C([0,T]; \R)$, $n_{\varepsilon} \geq C > 0$ for all values of $\varepsilon$ and for any $q \in \N_0$ there exists $C > 0$ such that 
	$$
	\sup_{t \in [0,T]} |n_{\varepsilon}(t)| \leq C \varepsilon^{q}
	$$
	for all values of $\epsilon$;
	
	\item[(b)] $n_{j,1,\varepsilon}, n_{0,\varepsilon} \in C([0,T];\C)$ and for any $q \in \N_0$ there exists $C > 0$ such that 
	$$
	\sum_{j=1}^{n} \sup_{t \in [0,T]} |n_{j,1,\varepsilon}(t)| + \sup_{t \in [0,T]} |n_{0,\varepsilon}(t)| \leq C \varepsilon^{q}
	$$
	for all values of $\epsilon$;
	
	\item[(c)] $l_{j,1,\varepsilon} \in \mathcal{B}^{\infty}(\R^{n})$, $|l_{j,1,\varepsilon}(x)| \leq C \langle x \rangle^{-2}$ and for any $q \in \N_0$ and any $\beta \in \N_0^n$ there exists $C > 0$ such that
	\begin{equation}\label{cnumber}
	\sup_{x \in \R^n} |\partial^{\beta}_{x} l_{j,1,\varepsilon}(x)| \leq C \varepsilon^{q}
	\end{equation}
	for all values  of $\varepsilon$;
	
	\item[(d)] $l_{0,\varepsilon} \in \mathcal{B}^{\infty}(\R^{n})$ and for any $q \in \N_0$ and any $\beta \in \N_0^n$ there exists $C > 0$ such that
	$$
	\sup_{x \in \R^n} |\partial^{\beta}_{x} l_{0,\varepsilon}(x)| \leq C \varepsilon^{q}
	$$
	for all values  of $\varepsilon$.
\end{itemize}
Let the net $(u'_{\varepsilon})_\varepsilon$ satisfies
\[
\begin{cases}
	S'_{\varepsilon} u'_\varepsilon(t,x) = f_\varepsilon(t,x) +p_{\varepsilon}(t,x), \quad t \in [0,T],\, x \in \R^n, \\
	u'_{\varepsilon}(0,x) = g_\varepsilon(x) + q_{\varepsilon}(x) \quad\qquad \,\,\,\,\, x \in \R^n,
\end{cases}
\]
where $(p_\varepsilon)_{\varepsilon}$ and $(q_\varepsilon)_{\varepsilon}$ are $H^{\infty}$-negligible nets. We now want to compare the two very weak solutions of $H^\infty$-type we obtained. 
We have
\[
\begin{cases}
	S_\varepsilon(u_\varepsilon-u'_\varepsilon)(t,x) = -p_{\varepsilon}(t,x)-(S_{\varepsilon}-S'_\varepsilon) u'_\varepsilon(t,x), \quad t \in [0,T],\, x \in \R^n, \\
	(u_\varepsilon-u'_\varepsilon)(0,x) =-q_{\varepsilon}(x), \qquad \qquad  \qquad \qquad \qquad \quad x \in \R^n.
\end{cases}
\]
By applying estimate \eqref{national_acrobat} with a suitable regularisation scale, we have that for all $m\ge 0$ there exist $C>0$ and $N\in\N_0$ such that 
\begin{equation}
	\label{est_neg_uniq}
	\|u_{\varepsilon}-u'_\varepsilon\|^{2}_{H^m} \leq  C\varepsilon^{-N} 
	\left\{ \|q_{\varepsilon}\|_{H^{m}} +  \int_{0}^{t} \|p_{\varepsilon}(\tau)+(S_{\varepsilon}-S'_\varepsilon) u'_\varepsilon(\tau)\|^{2}_{H^{m}} d\tau\right\}.
\end{equation}
Since $(p_\varepsilon)$ and $(q_\varepsilon)$ are both $H^{\infty}$-negligible, the coefficients of the operator $(S_\varepsilon-S'_\varepsilon)_\varepsilon$ satisfy negligible estimates and $(u'_{\varepsilon})_{\varepsilon}$ satisfies $H^{\infty}-$moderate estimates we conclude that the right-hand side of \eqref{est_neg_uniq} can be estimated by any positive power of $\varepsilon$. Hence, $(u_{\varepsilon}-u'_{\varepsilon})_{\varepsilon}$ is $H^\infty$-negligible. Summarising, we have proved the following uniqueness result.
\begin{Prop}\label{AFI}
	The very weak solutions of $H^\infty$-type of the Cauchy problem \eqref{CP} are unique in the following sense: negligible perturbations satisfying $(a), (b), (c)$ and $(d)$ on the regularisations of the equation coefficients and $H^{\infty}$-negligible perturbations of the regularisations of the initial data give $H^{\infty}$-negligible perturbation of the corresponding very weak solution.
\end{Prop}

\begin{Rem}
	Let us interpret the proposition above in the following particular case:
	\begin{itemize}
		\item $a \in C^{\infty}([0,T])$, $a$ is real-valued and 
		$$
		a(t) \geq C > 0;
		$$ 
		
		\item $b_j \in \mathcal{B}^{\infty}(\R^{n})$, $j = 0, 1, \ldots, n$ and 
		$$
		|b_j(x)| \leq C \langle x \rangle^{-2}, \quad j = 1, \ldots, n;
		$$
		
		\item $a_j \in C^{\infty}([0,T])$, $j = 0, 1,\dots, n$.
	\end{itemize}
	In this situation Proposition \ref{AFI} implies that the very weak solution does not depend on the mollifiers that we use to regularise the coefficients, provided that the mollifiers have all vanishing moments. Indeed, for $\ell = 1, 2$ consider 
	$$
	\rho^{(\ell)} \in \mathscr{S}(\R^n), \quad \int\rho^{(\ell)} = 1, \quad \int x^{\alpha}\rho^{(\ell)} = 0 \,\, \text{for any} \,\, \alpha \neq 0.
	$$
	Then, for $j = 1, \ldots, n$ define 
	$$
	b^{(\ell)}_{j,\varepsilon} = b_{j} \ast \rho^{(\ell)}_\varepsilon.
	$$
	From Proposition \ref{the_evil_that_men_do} we conclude that $(b^{(1)}_{j,\varepsilon} - b^{(2)}_{j,\varepsilon})_{\varepsilon}$ is $H^{\infty}-$negligible and it is easy to check that $|b^{(\ell)}_{j,\varepsilon}| \lesssim \langle x \rangle^{-2}$, in other words $(b^{(1)}_{j,\varepsilon} - b^{(2)}_{j,\varepsilon})_{\varepsilon}$ satisfies \eqref{cnumber}. Analogous considerations can be done for the regularisations of the coefficients $a$, $a_j$ and $b_0$ via mollifiers with all vanishing moments. Hence, Proposition \ref{AFI} gives that all very weak solutions obtained by all possible mollifiers with all vanishing moments are unique modulo $H^{\infty}$-negligible nets.
\end{Rem}

\subsection{Consistency with regular theory}

Suppose that the coefficients of the operator $S$ are regular as in the classical results, that is:
\begin{enumerate}
	\item[(i)] $a \in C([0,T]; \R)$ and $a(t) \geq C_{a} > 0$;
	
	\item[(ii)] $a_{j} \in C([0,T]; \C)$, $j =0, 1, \ldots, n$;
	
	\item[(iii)] $b_j \in \mathcal{B}^{\infty}(\R^{n})$, $j = 0, 1, \ldots, n$, and $|b_j(x)| \leq C \langle x \rangle^{-2}$, $j = 1, \ldots, n$.
\end{enumerate}
In this situation the Cauchy problem \eqref{CP} for the operator $S$ given in \eqref{true_operator_S} with data $f \in C([0,T];H^{\infty}(\R^{n}))$ and $g \in H^{\infty}(\R^{n})$ admits a unique solution $u \in C([0,T]; H^{\infty}(\R^{n}))$ satisfying an energy estimate like \eqref{national_acrobat} with a constant independent from $\varepsilon$. The goal here is to verify that the net obtained in Theorem \ref{mainn} converges to $u$ in the $H^{\infty}(\R^{n})$ topology. 

Let $u_{\varepsilon}$ be the solution of the Cauchy problem \eqref{regularized_CP}. Then $u - u_{\varepsilon}$ solves the following Cauchy problem
\begin{equation}\label{reise_reise}
	\begin{cases}
		S\{u-u_{\varepsilon}\}(t,x) = (f-f_\varepsilon)(t,x) + Q_{\varepsilon}u_{\varepsilon}(t,x), \quad t \in [0,T],\, x \in \R^{n}, \\
		\{u-u_{\varepsilon}\}(0,x) = (g-g_\varepsilon)(x), \qquad \qquad \qquad \quad \,\,\,\, x \in \R^{n},
	\end{cases}
\end{equation}
where 
$$
Q_{\varepsilon} = (a_{\varepsilon}(t)-a(t))\sum_{j=1}^{n}D^{2}_{x_j} + \sum_{j=1}^{n} (a_{j,\varepsilon}(t)b_{j,\varepsilon}(x) - a_j(t)b_{j}(x))D_{x_j} + (a_{0,\varepsilon}(t)b_{0,\varepsilon}(x)-a_0(t)b_0(x)).
$$

So, we have the following estimate
\begin{align*}\label{nebel}
	\|&\{u-u_{\varepsilon}\}(t)\|^{2}_{H^m} \leq  C_{m} \left\{ \|g-g_{\varepsilon}\|^{2}_{H^m} + \int_{0}^{t} \|(f-f_\varepsilon)(\tau) + Q_{\varepsilon} u_{\varepsilon}(\tau)\|^{2}_{H^{m}} d\tau  \right\} \\
	&\leq C_{m}\|g-g_{\varepsilon}\|^{2}_{H^m} + C_m \int_{0}^{t} \|(f-f_\varepsilon)(\tau)\|^{2}_{H^m}  d\tau + C_m \int_{0}^{t} |a_\varepsilon(\tau) - a(\tau)|^{2} \|u_{\varepsilon}(\tau)\|^{2}_{m+2} d\tau \\
	&+ C_m \int_{0}^{t} \sum_{j=1}^{n} \max_{|\alpha| \leq m} \sup_{x \in \R^{n}}|a_{j,\varepsilon}(\tau)D^{\alpha}_{x}b_{j,\varepsilon}(x) - a_j(\tau)D^{\alpha}_{x}b_{j}(x)|^{2} \|u_{\varepsilon}(\tau)\|^{2}_{m+1} d\tau \\
	&+ C_m \int^{t}_{0} \max_{|\alpha| \leq m}\sup_{x\in\R^n}|a_{0,\varepsilon}(\tau)D^{\alpha}_{x}b_{0,\varepsilon}(x)-a_0(\tau)D^{\alpha}_{x}b_0(x)|^{2} \|u_{\varepsilon}(\tau)\|^{2}_{m} d\tau.
\end{align*}
Now we observe the following:
\begin{itemize}
	\item[(i)] Since in this case we are regularising \textit{regular} functions, it is easy to conclude that the regularisations $a_\varepsilon$, $a_{j,\varepsilon}$ and $b_{j,\varepsilon}$ satisfy uniform estimates with respect to $\varepsilon$. So, in this particular case, we obtain estimate \eqref{national_acrobat} uniformly with respect to $\varepsilon$. Hence, for any $m \in \N_0$, 
	$$
	\|u_{\varepsilon}(\tau)\|^{2}_{m} \leq C_{m}, \quad \forall \varepsilon \in (0,\varepsilon_0).
	$$
	
	\item[(ii)] We have the following uniform convergence
	$$
	\sup_{t \in [0,T]} |a_{\varepsilon}(t) - a(t)| \to 0 \quad \text{and} \quad \sup_{t \in [0,T]} |a_{j,\varepsilon}(t) - a_{j}(t)| \to 0.
	$$
	
	\item[(iii)] We also have $b_{j,\varepsilon} \to b_j$ in $\mathcal{B}^{\infty}(\R^{n})$ as $\varepsilon \to 0$, i.e., for any $\beta \in \N^n_0$ it holds 
	$$
	\sup_{x \in \R^{n}}|\partial^{\beta}_{x}b_{j,\varepsilon}(x) - \partial^{\beta}_{x}b_{j}(x)| \to 0, \quad \text{as} \, \varepsilon \to 0.
	$$
	
	\item[(iv)] $(g-g_\varepsilon)_{\varepsilon}$ and $((f-f_{\varepsilon})(t))_{\varepsilon}$ are $H^{\infty}-$negligible nets, provided that we regularise the initial data $g$ and $f$ via mollifiers with all vanishing moments.
\end{itemize}
We therefore conclude $u_{\varepsilon} \to u$ in $C([0,T];H^{\infty}(\R^{n}))$. 

\begin{Prop}
	Assume that $a \in C([0,T];\R)$, $a$ never vanishes, $a_j \in C([0,T];\C)$, $b_j \in \mathcal{B}^{\infty}(\R^{n})$ for $j = 0, 1, \ldots, n$ and $|b_{j}(x)| \leq C \langle x \rangle^{-2}$ for $j = 1, \ldots, n$. Assume moreover $g \in H^{\infty}(\R^n)$ and $f \in C([0,T];H^{\infty}(\R^n))$. Then, if we use mollifiers with all vanishing moments to regularise the Cauchy data, any very weak solution converges to the unique classical solution in the space $C([0,T];H^{\infty}(\R^{n}))$. In particular, in the classical case, the limit of very weak solutions always exist and does not depend on the regularisation. 
\end{Prop}

\end{document}